%% file: paper.tex
\documentclass[10pt]{article}

\usepackage{cite}

\usepackage{amsfonts, amssymb, amsmath, amsthm, epsf}

\title{Symmetric periodic solutions of parabolic problems  with hysteresis}
\author{Pavel Gurevich\thanks{partially supported by  DFG project
SFB 555 and the RFBR project 10-01-00395-a}  \\
Institute for Mathematics I, Free University of Berlin, \\
Arnimallee 3, 14195 Berlin, Germany. \\ gurevichp@gmail.com
\and Sergey Tikhomirov\thanks{partially supported by NSC (Taiwan)  97-2115-M-002 -011 -MY2 and CNPq (Brazil)}  \\
Departamento de Matematica\\
Pontificia Universidade Catolica do Rio de Janeiro\\
Rua Marques de Sao Vicente, 225\\
Edificio Cardeal Leme, sala 862\\
Gavea - Rio de Janeiro - Brazil CEP 22453-900\\
sergey.tikhomirov@gmail.com}
\date{\today}

\theoremstyle{plain}
\newtheorem{theorem}{Theorem}[section]
\newtheorem{lemma}{Lemma}[section]
\newtheorem{corollary}{Corollary}[section]
\newtheorem{condition}{Condition}[section]

\theoremstyle{definition}
\newtheorem{definition}{Definition}[section]
\newtheorem{example}{Example}[section]
\newtheorem{remark}{Remark}[section]

\numberwithin{equation}{section} \numberwithin{figure}{section}

\renewcommand{\Im}{{\rm Im\,}}
\renewcommand{\Re}{{\rm Re\,}}

\newcommand{\mes}{{\rm mes\,}}

\newcommand{\Span}{{\rm Span}}

\newcommand{\const}{{\rm const}}
\newcommand{\dist}{{\rm dist}}
\renewcommand{\phi}{{\varphi}}

\newcommand{\tr}{{\rm tr\,}}

\newcommand{\cV}{{\mathcal V}}

\newcommand{\cH}{{\mathcal H}}

\newcommand{\cW}{{\mathcal W}}

\newcommand{\cU}{{\mathcal U}}

\newcommand{\bA}{{\mathbf A}}
\newcommand{\bE}{{\mathbf E}}

\newcommand{\bP}{{\mathbf P}}
\newcommand{\bR}{{\mathbf R}}

\newcommand{\bB}{{\mathbf B}}

\newcommand{\bI}{{\mathbf I}}

\newcommand{\bM}{{\mathbf M}}

\newcommand{\bv}{{\mathbf v}}
\newcommand{\bw}{{\mathbf w}}

\newcommand{\bt}{{\mathbf t}}
\newcommand{\bz}{{\mathbf z}}
\newcommand{\bphi}{{\boldsymbol{\varphi}}}
\newcommand{\bpsi}{{\boldsymbol{\psi}}}

\newcommand{\bbR}{{\mathbb R}}
\newcommand{\bbN}{{\mathbb N}}

\newcommand{\bbJ}{{\mathbb J}}

\let\phi=\varphi

\begin{document}

\maketitle

\begin{abstract}
We consider the heat equation in a multidimensional domain with
nonlocal hysteresis feedback control in a boundary condition.
Thermostat is our prototype model. We construct all periodic
solutions with exactly two switching on the period and study their
stability. Coexistence of several periodic solutions with
different stability properties is proved to be possible. A
mechanism of appearance and disappearance of periodic solutions is
investigated.
\end{abstract}

\setcounter{tocdepth}{2} \tableofcontents

\input intro.tex
\input sect2.tex
\input sect3.tex
\input sect4.tex

\bigskip

\input bibl.tex
\end{document}

%% file: intro.tex
\section{Introduction}

Hysteresis operators  arise in mathematical description of various
physical processes~\cite{KrasnBook,Visintin,BrokateSprekels}.
Models with hysteresis for  ordinary differential equations  were
considered by many authors (see e.g.,~\cite{Alt, KrasnBook,
Seidman, Feckan,VarigondaGeorgiou,
BlimanKrasnos,MackiNistriZecca}). Partial differential equations
with hysteresis have also been actively studied during the last
decades (see~\cite{Visintin,BrokateSprekels} and the references
therein). The primary focus has been on the well-posedness of the
corresponding problems and related issues (existence of solutions,
uniqueness, regularity, etc.). However, many questions remain
open, especially those related to the periodicity and long-time
behavior of solutions.

In this paper, we deal with parabolic problems containing a
discontinuous hysteresis operator in the boundary condition. Such
problems describe  processes of thermal control arising in
chemical reactors and climate control systems. The temperature
regulation in a domain is performed via heating (or cooling)
elements on the boundary of the domain. The regime of the heating
elements on the boundary is based on the registration of thermal
sensors inside the domain and obeys a hysteresis law.

Let $v(x,t)$ denote  the temperature at the point $x$ of a bounded
domain $Q\subset\bbR^n$ at the moment $t$. We define the {\it mean
temperature} $\hat v(t)$ by the formula
$$
\hat v(t)=\int_Q m(x)v(x,t)\,dx,
$$
where $m$ is a given function from the Sobolev space $H^1(Q)$ (see
Condition~\ref{condKm} for another technical assumption on
$m(x)$).

In our prototype model, we assume that the function $v(x,t)$
satisfies the heat equation
\begin{equation}\label{eqHeatEqu}
v_t(x,t) =\Delta v(x,t) \quad (x\in Q,\ t>0)
\end{equation}
and a boundary condition which involves a hysteresis operator
$\cH$ depending on the mean temperature $\hat v$.

 The hysteresis $\cH(\hat v)(t)$ is
defined as follows (cf.~\cite{KrasnBook,Visintin} and the accurate
definition and Fig.~\ref{figHyst} in Sec.~\ref{secSetting}). One
fixes two temperature thresholds $\alpha$ and~$\beta$
($\alpha<\beta$). If $\hat v(t)\le\alpha$, then $\cH(\hat v)(t)=1$
(the heating is switched on); if $\hat v(t)\ge\beta$, then
$\cH(\hat v)(t)=-1 $ (the cooling is switched on); if the mean
temperature $\hat v(t)$ is between $\alpha$ and~$\beta$, then
$\cH(\hat v)(t)$ takes the same value as ``just before.'' We say
that the hysteresis operator {\it switches} when it jumps from $1$
to $-1$ or from $-1$ to $1$. The corresponding time moment is
called the {\it switching moment.} Note that the hysteresis
phenomenon takes place along with the nonlocal effect caused by
averaging of the function $v(x,t)$ over $Q$.

To be definite, let us assume that one regulates the heat flux
through the boundary $ \partial Q$. Then the boundary condition is
of the form
\begin{equation}\label{eqBoundaryCondition}
\frac{\partial v}{\partial\nu} = K(x) \cH(\hat v)(t)   \quad
(x\in\partial Q,\ t>0),
\end{equation}
where $\nu$ is the outward normal to $\partial Q$ at the point
$x$, $K$ is a given smooth real-valued function (distribution of
the heating elements on the boundary).

A similar mathematical model   was originally proposed
in~\cite{GlSpr-SIAM,GlSpr-JIntEqu}. Generalizations to various
phase-transition problems with hysteresis    were studied
in~\cite{CGS,FH,HNS,Krejci, BrokateSprekels}. Some related issues
of optimal control were considered in~\cite{BrokateFriedman}. The
most important questions here concern the existence and uniqueness
of solutions, the existence of periodic solutions, and long-time
behavior of solutions. The latter two questions are especially
difficult.

In the case of a {\it one-dimensional} domain $Q$ (a finite
interval, $n=1$),  the periodicity was studied
in~\cite{FriedmanJiang-CPDE,Pruess,GoetzHoffmannMeirmanov,Kopfova}.
Problems with hysteresis on the boundary of a {\it
multidimensional} domain ($n\ge 2$) turn out to be much more
complicated. Although one can relatively easily prove the
existence (and sometimes uniqueness) of solutions, the issue of
finding periodic solutions is still an open question. The main
difficulty here is related to the fact that, in general, the
solution does not depend on the initial data continuously. The
reason is that the solution may intersect the ``switching''
hyperplane $\{\hat \phi=\alpha\}$ or $\{\hat \phi=\beta\}$
nontransversally (cf.~\cite{BlimanKrasnos,VarigondaGeorgiou},
where the same phenomenon occurs for ordinary differential
equations). This leads to  discontinuity of the corresponding
Poincar\'e map. As a result, most methods based on fixed-point
theorems do not apply to the  Poincar\'e map.

One possible way to overcome the nontransversality is to consider
a continuous model of the hysteresis operator. This was done
in~\cite{GurJaegerPreisach}, where a thermocontrol problem with
the Preisach hysteresis operator in the boundary condition was
considered and the existence of periodic solutions and global
attractors were established. Note that the periodicity and the
long-time behavior of solutions were also studied   in~\cite{Xu,
KenmochiVisintin} in the situation where a  hysteresis operator
enters a parabolic equation itself  (see
also~\cite{Visintin_Annali} and the references therein).

The first results about periodic solutions of thermocontrol
problems in {\it multidimensional} domains with {\it
discontinuous} hysteresis were obtained
in~\cite{GurJaegerSkubSIAM}. In~\cite{GurHystUncond},  a new
approach  was proposed. It is based on regarding the problem as an
infinite-dimensional dynamical system. By using the Fourier
method, one can reduce the boundary-value problem for the
parabolic equation to infinitely many ordinary differential
equations, whose solutions are coupled with each other via the
hysteresis operator.

In~\cite{GurHystUncond}, the existence of a unique periodic
solution of the thermocontrol problem is proved for sufficiently
large $\beta-\alpha$. This periodic solution possesses certain
symmetry, is stable, and is a global attractor. A similar result
was established for arbitrary $\alpha$ and $\beta$, but $m(x)$
being close to a constant. The idea  was to find an invariant
region for the corresponding Poincar\'e map and prove that the
Poincar\'e map is continuous on that region. This turns out to be
true for sufficiently large $\beta-\alpha$. However, one can
construct examples with small $\beta-\alpha$, where an invariant
region exists and even is  an attracting set, but the Poincar\'e
map is not continuous on it.

In the present paper, we will show that the requirement for
$\beta-\alpha$ to be large enough is essential. We will prove that
if $\beta-\alpha$ is small, then unstable periodic solutions may
appear. In particular, they may have a saddle structure. To
construct those solutions, we will develop a general procedure
which  yields {\it all} periodic solutions (with  two switchings
on the period) in an explicit form. This procedure works even in
the presence of discontinuity caused by the above
nontransversality. In particular, it allows one to find periodic
solutions on which the Poincar\'e map is discontinuous.

To study   stability of periodic solutions (in particular, to find
unstable ones), we propose a method which allows one to reduce the
original system to an invariant subsystem. The dimension of this
subsystem is equal to the number of  nonvanishing modes in the
Fourier decomposition of the $m(x)$. If  $m(x)$ has finitely many
nonvanishing modes, then  one can explicitly write down the
linearization of the reduced system and find all the eigenvalues.
They provide complete information about the stability of the
 periodic solution.

The invariant subsystem corresponding to the nonvanishing modes of
$m(x)$  is called {\it guiding}. The remaining subsystem is called
{\it guided}. We prove that the full system (i.e., the original
problem) has a periodic solution whenever the guiding system has
one. Moreover, the periodic solution of the full system is a
global attractor (is stable, uniformly exponentially stable)
whenever the periodic solution of the guiding system possesses
those properties. We call these results {\it conditional
existence} of periodic solutions, {\it conditional attractivity},
and {\it conditional stability}, respectively. The above
``guiding-guided'' decomposition is a result of independent
interest. It generalizes the results of~\cite{GurJaegerSkubSIAM},
where $m(x)\equiv\const$ (in our terminology, this corresponds to
$m(x)$ which has only one nonvanishing mode).

The paper is organized as follows. In Sec.~\ref{secSetting}, we
define the hysteresis operator, formulate the problem,  introduce
a notion of  solution,   recall some properties of the solutions,
and reduce the problem to an infinite-dimensional dynamical
system. In the end of Sec.~\ref{secSetting}, we define the guiding
and the guided subsystems and introduce the corresponding
decomposition of the phase space (the Sobolev space $H^1(Q)$).
Most results of this section are proved in~\cite{GurHystUncond}.

In Sec.~\ref{secPeriodic}, we give a notion of periodic solution
with two switchings on the period. By using the Poincar\'e maps of
the guiding system and the full system, we prove conditional
existence of periodic solutions, conditional attractivity, and
conditional stability. The latter two results are proved under
assumption that the periodic solution of the guiding system
intersects the hyperplanes $\{\hat\phi=\alpha\}$ and
$\{\hat\phi=\beta\}$ at the switching moments transversally. The
transversality implies the continuity (and even the Fr\'echet
differentiability) of the Poincar\'e maps in a  neighborhood of
the periodic solution. However, we require neither that this
neighborhood be invariant under the Poincar\'e map, nor that the
Poincar\'e map be continuous in a (bigger) invariant neighborhood
(which exists due to~\cite{GurHystUncond}).

In Sec.~\ref{secPeriodicReduction}, we show that any periodic
solution with two switchings on the period possesses a symmetry in
the phase space. By using this symmetry, we develop  an algorithm
which allows us
\begin{enumerate}

\item  to construct {\it all} periodic solutions (with two
switchings on the period) in an explicit form for any given
$\alpha$ and $\beta$;

\item to find a sufficient condition under which periodic
solutions exist for all sufficiently small $\beta-\alpha$.

\item to define bifurcation points where periodic solutions may
appear or disappear; a role of a bifurcation parameter is played
either by the period or by the difference $\beta-\alpha$;

\end{enumerate}

 Furthermore, using the
results about the guiding-guided decomposition from
Sec.~\ref{secPeriodic}, we construct examples in which  periodic
solutions are stable or unstable, respectively. In the
``unstable'' case, we show that they may have a saddle structure.

As a conclusion, we note that the developed method  can also be
applied to the study of the Dirichlet or Robin boundary
conditions. Moreover, one can study the problem where the heat
flux through the boundary (in the case of the Neumann boundary
condition) changes continuously. Mathematically, this means that
the boundary condition~\eqref{eqBoundaryCondition} is replaced by
$$
\frac{\partial v}{\partial\nu} = K(x) u(t)   \quad (x\in\partial Q
,\ t>0),
$$
$$
au'(t)+u(t)=\cH(\hat v)(t)
$$
with $a>0$
(cf.~\cite{GlSpr-SIAM,GlSpr-JIntEqu,Pruess,GurJaegerSkubSIAM}).

%% file: sect2.tex
\section{Setting of the Problem. Reduction to Infinite Dynamical System}\label{secSetting}

\subsection{Setting of the problem}\label{subsecSetting} Let $Q\subset\bbR^n$ ($n\ge 1$) be a bounded domain
with smooth boundary. Let $L_2=L_2(Q)$. Denote by $H^1=H^1(Q) $
  the Sobolev space with the norm
$$
\|\psi\|_{H^1}=\left(\int_Q
(|\psi(x)|^2+|\nabla\psi(x)|^2)\,dx\right)^{1/2}.
$$
Let $H^{1/2}=H^{1/2}(\partial Q)$  be the space of traces on
$\partial Q$ of the functions from $H^1$.

Consider the sets $Q_T=Q\times(0,T)$ and $\Gamma_T=\partial Q
\times(0,T)$, $T>0$. Fix functions $K\in H^{1/2}$ and $m\in H^1$
and real numbers $\alpha$ and $\beta$, $\beta>\alpha$.

For any function $\phi(x)$ or $v(x,t)$ ($x\in Q$, $t\ge0$), the
symbol $\hat{\phantom{a}}$ will refer to the ``average'' of the
function:
$$
\hat \phi=\int_Q m(x) \phi(x)\,dx,\qquad \hat v(t)=\int_Q m(x)
v(x,t)\,dx.
$$

Let $v(x,t)$ denote the temperature at the point $x\in Q$ at the
moment $t\ge0$ satisfying the heat equation
\begin{equation}\label{eq2.1}
v_t(x,t) =\Delta v(x,t) \quad ((x,t)\in Q_T)
\end{equation}
with the initial condition
\begin{equation}\label{eq2.2}
v|_{t=0}=\phi(x)\quad (x\in Q)
\end{equation}
and the boundary condition
\begin{equation}\label{eq2.3}
  \frac{\partial v}{\partial\nu}\Big|_{\Gamma_T} = K(x) \cH(\hat v)(t)  \quad
((x,t)\in\Gamma_T).
\end{equation}
Here  $\nu$ is the outward normal to $\Gamma_T$ at the point
$(x,t)$ and $\cH$ is a hysteresis operator, which we now define.

We denote by $BV(0,T)$ the Banach space of real-valued functions
having finite total variation on the segment $[0,T]$ and by
$C_r[0,T )$ the linear space of functions which are continuous on
the right in $[0,T )$. We  introduce the {\it hysteresis
operator\/} (cf.~\cite{KrasnBook,Visintin})
$$
\cH: C[0,T]\to BV(0,T)\cap C_r[0,T )
$$
by the following rule. For any $g\in C[0,T]$, the function
$h=\cH(g):[0,T]\to\{-1,1\}$ is defined as follows. Let
$X_t=\{t'\in(0,t]:g(t')=\alpha\ \text{or } \beta\}$; then
$$
h(0)=\begin{cases}
1 & \text{if } g(0 )< \beta,\\
-1 & \text{if } g(0)\ge \beta
\end{cases}
$$
and for $t\in(0,T]$
$$
h(t)=\begin{cases}
h(0) & \text{if } X_t=\varnothing,\\
1 & \text{if } X_t\ne\varnothing\ \text{and } g(\max X_t)=\alpha, \\
-1 & \text{if } X_t\ne\varnothing\ \text{and } g(\max X_t)=\beta
\end{cases}
$$
\begin{figure}[ht]
      {\ \hfill\epsfxsize55mm\epsfbox{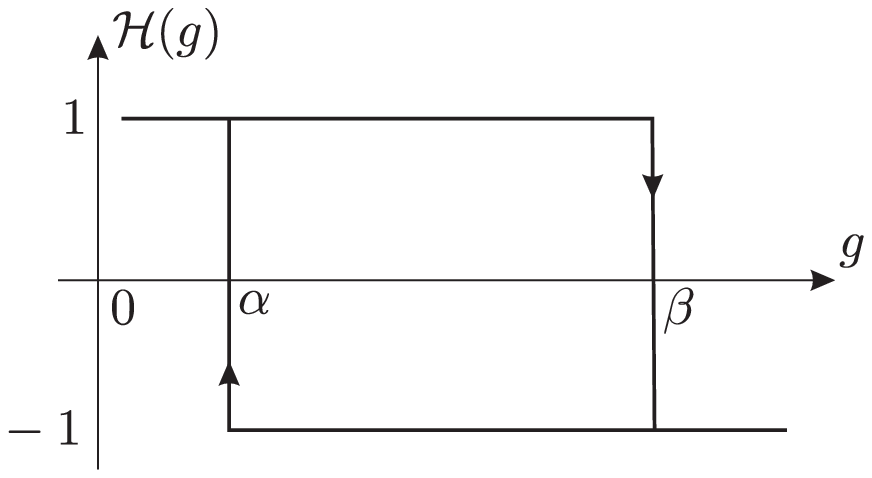}\hfill\ }
        \caption{The hysteresis operator $\cH$}
        \label{figHyst}
\end{figure}
(see Fig.~\ref{figHyst}). A point $\tau$ such that
$\cH(g)(\tau)\ne\cH(g)(\tau-0)$ is called a {\it switching moment
of $\cH(g)$}.

We assume throughout that the following condition holds.
\begin{condition}\label{condKm}
The coefficient $K(x)$ in the boundary condition~\eqref{eq2.3} and
the weight function $m(x)$ satisfy
\begin{equation}\label{eqKmnot0}
\int_{\partial Q} K(x)\,d\Gamma>0,\qquad \int_Q m(x)\,dx>0.
\end{equation}
\end{condition}

\begin{remark}
From the physical viewpoint, the function $K(x)$ characterizes the
density of the heating (or cooling) elements on the boundary and
$m(x)$ characterizes the density of thermal sensors in the domain.
Clearly, inequalities~\eqref{eqKmnot0} hold in the physically
relevant case  $K(x)\ge0$ for a.e. $x\in \partial Q$,
$K(x)\not\equiv 0$, and $m(x)\ge0$ for a.e. $x\in Q$,
$m(x)\not\equiv 0$.
\end{remark}

\subsection{Functional spaces and the solvability of the problem}

For any Banach space $B$, denote by $C([a,b];B)$ ($a<b$) the space
of $B$-valued functions continuous on the segment $[a,b]$ with the
norm
$$
\|u\|_{C([a,b];B)}=\max\limits_{t\in[a,b]}\|u(t)\|_{B}
$$
and by $L_2((a,b);B)$ the space of $L_2$-integrable $B$-valued
functions with the norm
$$
\|u\|_{L_2((a,b);B)}=\left(\int_a^b \|u(t)\|_B\,dt\right)^{1/2}.
$$

We introduce the anisotropic Sobolev space $
H^{2,1}(Q\times(a,b))=\{v\in L_2((a,b); H^2): v_t\in
L_2((a,b);L_2)\}$  with the norm
$$
\|v\|_{H^{2,1}(Q\times(a,b))} =
  \left( \int_a^b \|v(\cdot,t)\|_{H^2 }^2 \,dt+ \int_a^b\left\|
v_t(\cdot,t) \right\|_{L_2}^2 \,dt\right)^{1/2}.
$$

Taking into account the results of the interpolation theory (see,
e.g.,~\cite[Chap.~1, Secs.~1--3, 9]{LM}, we make the following
remarks.

\begin{remark}\label{remWC}
The continuous embedding $H^{2,1}(Q\times(a,b))\subset
C([a,b],H^1)$ takes place. Furthermore,  for any $v\in
H^{2,1}(Q\times(a,b))$ and $\tau\in[a,b]$, the trace
$v|_{t=\tau}\in H^1$ is well defined and is a bounded operator
from $H^{2,1}(Q\times(a,b))$ to $H^1$.
\end{remark}

\begin{remark}\label{remWConcat}
Consider two functions  $v_1\in H^{2,1}(Q\times(a,b))$ and $v_2\in
H^{2,1}(Q\times(b,c))$, where $a<b<c$. Let
$v(\cdot,t)=v_1(\cdot,t)$ for $t\in (a,b)$ and
$v(\cdot,t)=v_2(\cdot,t)$ for $t\in (b,c)$. Then $v\in
H^{2,1}(Q\times(a,c))$ if and only if $v_1|_{t=b}=v_2|_{t=b}$.
\end{remark}

\begin{definition}\label{defStrongHyst}
A   function $ v(x,t) $ is called a solution of
problem~\eqref{eq2.1}--\eqref{eq2.3}  in~$Q_T$  with the initial
data $\phi\in H^1$ if $v\in H^{2,1}(Q_T)$ and $v$ satisfies
Eq.~\eqref{eq2.1} a.e.  in~$Q_T$
    and conditions~\eqref{eq2.2}, \eqref{eq2.3} in the sense of traces.
\end{definition}
\begin{definition}
We say that  $v(x,t)$ ($t\ge0$) is a {\it solution  of
problem~\eqref{eq2.1}--\eqref{eq2.3} in~$Q_\infty$} if it is a
solution in $Q_T$ for all $T>0$.
\end{definition}

The following result about the solvability of
problem~\eqref{eq2.1}--\eqref{eq2.3} is proved in~\cite[Theorem
2.2]{GurHystUncond}.
\begin{theorem}\label{tExistenceUniqueness}
Let $\phi\in H^1$ and $\|\phi\|_{H^1}\le R$ {\rm (}$R>0$ is
arbitrary{\rm )}. Then there exists a unique solution $v$ of
problem~\eqref{eq2.1}--\eqref{eq2.3} in~$Q_\infty$ and   the
following holds for any $T>0$.
\begin{enumerate}
\item One has
\begin{equation}\label{eqExistenceUniqueness}
\|v(\cdot, t)\|_{H^1}\le c_0 \|v\|_{H^{2,1}(Q_T)}\le
c_1(\|\phi\|_{H^1}+\|K\|_{H^{1/2}}),
\end{equation}
where $c_0=c_0(T)>0$ and $c_1=c_1(T)>0$ do not depend on $\phi$
and $R${\rm ;} \item  The interval $(0,T]$ contains no more than
finitely many switching moments $t_1<t_2<{\dots}<t_J$
 of $\cH(\hat v)$. Moreover,
\begin{equation}\label{eqExistenceUniquenessMild'}
t_i-t_{i-1}\le t^*+\dfrac{2(\beta-\alpha)}{m_0K_0},\quad
i=1,2,\dots,
\end{equation}
where $t^*$ depends on $m$ and $R$ but does not depend on
$\phi,T,\alpha,\beta${\rm ;}
\begin{equation}\label{eqt*}
t_i-t_{i-1}\ge \tau^*,\quad i=
\begin{cases}
 1,2,\dots,J\quad \text{if }\hat\phi\le\alpha\ \text{or }\hat\phi\ge\beta,\\
 2,3,\dots,J\quad \text{if }\alpha<\hat\phi<\beta,
\end{cases}
\end{equation}
where
\begin{equation}\label{eqtau*ge}
\tau^*= \const\dfrac{(\beta-\alpha)^2}{\|m\|_{L_2}^2}
\end{equation}
 with $\const>0$ depending on $R$ rather than on $m,\phi,T,\alpha,\beta${\rm ;}
 in~\eqref{eqExistenceUniquenessMild'} and~\eqref{eqtau*ge},  $t_0=0$.

\end{enumerate}
\end{theorem}

\subsection{Reduction to   infinite-dimensional dynamical system}

Due to Theorem~\ref{tExistenceUniqueness}, the study of the
solutions of problem~\eqref{eq2.1}--\eqref{eq2.3} with hysteresis
can be reduced to the study of the solutions of parabolic problems
without  hysteresis by considering the time intervals between the
switching moments $t_i$.

Thus, if $\cH(\hat v)(t)\equiv 1$, then
problem~\eqref{eq2.1}--\eqref{eq2.3} takes the form
\begin{equation}\label{eqExistenceUniquenessNoHyst1}
v_t(x,t) =\Delta v(x,t) \quad ((x,t)\in Q_T),
\end{equation}
\begin{equation}\label{eqExistenceUniquenessNoHyst2}
v(x,0)=\phi(x)\quad (x\in Q),
\end{equation}
\begin{equation}\label{eqExistenceUniquenessNoHyst3}
  \frac{\partial v}{\partial\nu}\Big|_{\Gamma_T} =   K(x)    \quad
((x,t)\in\Gamma_T).
\end{equation}
If $\cH(\hat v)(t)\equiv -1$, one should replace $K(x)$ by $-K(x)$
in~\eqref{eqExistenceUniquenessNoHyst3}.

\begin{definition} A  function $ v(x,t) $ is called a    solution of
problem~\eqref{eqExistenceUniquenessNoHyst1}--\eqref{eqExistenceUniquenessNoHyst3}
in~$Q_T$ if $v\in H^{2,1}(Q_T)$ and $v$ satisfies
Eq.~\eqref{eqExistenceUniquenessNoHyst1} a.e. in~$Q_T$ and
conditions~\eqref{eqExistenceUniquenessNoHyst2},
\eqref{eqExistenceUniquenessNoHyst3} in the sense of traces.
\end{definition}

It is well known that there is a unique solution $v\in
H^{2,1}(Q_T)$ of
problem~\eqref{eqExistenceUniquenessNoHyst1}--\eqref{eqExistenceUniquenessNoHyst3}.

Now we give a convenient representation of solutions of
problem~\eqref{eqExistenceUniquenessNoHyst1}--\eqref{eqExistenceUniquenessNoHyst3}
in terms of the Fourier series with respect to the eigenfunctions
of the Laplacian.

Let $\{\lambda_j\}_{j=0}^\infty$ and $\{e_j(x)\}_{j=0}^\infty$
denote the sequence of eigenvalues and the corresponding system of
real-valued eigenfunctions (infinitely differentiable in
$\overline Q$) of the spectral problem
\begin{equation}\label{eqEigenProblem}
-\Delta e_j(x)=\lambda_j e_j(x)\quad (x\in Q),\qquad
\frac{\partial e_j}{\partial\nu}\Big|_{\partial Q} = 0.
\end{equation}

It is well known that $0=\lambda_0<\lambda_1\le\lambda_2\le{\dots}
\le\lambda_j\le\dots$, $e_0(x)\equiv (\mes Q)^{-1/2}>0$, and the
system of eigenfunctions $\{e_j\}_{j=0}^\infty$ can be chosen to
form an orthonormal basis for $L_2$. Then,  the functions
$e_j/\sqrt{\lambda_j+1}$ form an orthonormal basis for $H^1$.

\begin{remark}\label{remAsymp}
In what follows, we will use the well-known asymptotics for the
eigenvalues $\lambda_j= L j^{2/n}+o(j^{2/n})$ as $j\to+\infty$
($L>0$ and $n$ is the dimension of $Q$).
\end{remark}

Any function $\psi\in L_2$  can be expanded into the Fourier
series with respect to $e_j(x)$, which converges in $L_2$:
\begin{equation}\label{eqFourierPsiL2}
\psi(x) =\sum\limits_{j=0}^\infty \psi_{j} e_j(x),\qquad
\|\psi\|_{L_2}^2=\sum\limits_{j=0}^\infty |\psi_j|^2,
\end{equation}
where $ \psi_{j}=\int_Q \psi(x) e_j(x)\,dx.$ If $\psi\in H^1$,
then the first series in~\eqref{eqFourierPsiL2}   converges to
$\psi$ in $H^1$ and
\begin{equation}\label{eqFourierPsiH1}
\|\psi\|_{H^1}^2=\sum\limits_{j=0}^\infty (1+\lambda_j)
|\psi_j|^2.
\end{equation}

Denote
\begin{equation}\label{eqmjKj}
\begin{gathered}
m_j=\int_Q m(x)e_j(x)\,dx,\qquad K_j=\int_{\partial Q}
K(x)e_j(x)\,dx \qquad (j=0,1,2,\dots).
\end{gathered}
\end{equation}
Note that $m_0,K_0>0$ due to Condition~\ref{condKm}. We also note
that $K_j$ are not the Fourier coefficients of $K(x)$. However,
the following is proved  in~\cite{GurHystUncond}:
\begin{equation}\label{eqKjConvergence}
\sum\limits_{j=1}^\infty\left(\dfrac{|K_j|^2}{\lambda_j^2}+\dfrac{|K_j|^2}{\lambda_j}\right)\le
c\|K\|_{H^{1/2}}^2,
\end{equation}
where $c>0$ does not depend on $K$.

\begin{remark}\label{remmjKj}
Using~\eqref{eqFourierPsiH1} and~\eqref{eqKjConvergence}, we
obtain the estimate
$$
\sum\limits_{j=1}^\infty|m_jK_j|\le c\|m\|_{H^1}\|K\|_{H^{1/2}},
$$
which will often be used later on.
\end{remark}

The numbers $m_j$ and $K_j$ play an essential role when one
describes the thermocontrol problem in terms of an
infinite-dimensional dynamical system. The following result is
true (see~\cite[Lemma~2.2]{GurHystUncond}).

\begin{theorem}\label{thExistenceUniquenessNoHystMild}
Let $\phi\in H^1$. Then the following assertions hold.
\begin{enumerate}
\item
The  solution $v$ of
problem~\eqref{eqExistenceUniquenessNoHyst1}--\eqref{eqExistenceUniquenessNoHyst3}
can be represented as the series
\begin{equation}\label{eqvFourier}
v(x,t)=\sum\limits_{j=0}^\infty v_j(t)e_j(x),\quad t\ge0,
\end{equation}
where $v_j(t)=\int_Q v(x,t)e_j(x)\,dx$ and $v_j(t)$ satisfy the
Cauchy problem
\begin{equation}\label{eqExistenceUniquenessNoHyst5''}
\dot v_j(t)=-\lambda_j  v_j(t)+K_j,\qquad
 v_j(0)=\phi_j\qquad (\dot{\phantom{v}}=d/dt,\ j=0,1,2\dots).
\end{equation}
 The series in~\eqref{eqvFourier} converges in  $H^1$ for all
$t\ge0$.
\item The mean temperature $\hat v(t)$ is represented by the absolutely convergent
series
\begin{equation}\label{eqMeanTempFourier}
\hat v(t)=\sum\limits_{j=0}^\infty m_jv_j(t),\quad t\ge0,
\end{equation}
which is   continuously differentiable for $t>0$.
\end{enumerate}
\end{theorem}

\begin{remark}\label{rExplicitVj}
In what follows, we will also use the explicit formulas for the
solutions of Eqs.~\eqref{eqExistenceUniquenessNoHyst5''}
$$
v_0(t)=\phi_0+K_0t,\qquad
v_j(t)=\left(\phi_j-\dfrac{K_j}{\lambda_j}\right)e^{-\lambda_j
t}+\dfrac{K_j}{\lambda_j},\quad j=1,2,\dots.
$$
\end{remark}

Formally, relations~\eqref{eqExistenceUniquenessNoHyst5''} can be
obtained by multiplying~\eqref{eqExistenceUniquenessNoHyst1} by
$e_j(x)$, integrating by parts over $Q$, and substituting
$v(x,t)=\sum\limits_{j=0}^\infty v_j(t)e_j(x)$.  The rigorous
proof is given in~\cite{GurHystUncond}.

A geometrical interpretation of the dynamics of $v_0(t),
v_1(t),\dots$ is as follows. We choose the orthonormal basis in
$L_2$ (which is orthogonal in $H^1$) consisting of the
eigenfunctions $e_0,e_1,e_2,\dots$. Then, in the coordinate form,
we have
$$
e_0=(1,0,0,0,\dots),\quad e_1=(0,1,0,0,\dots),\quad
e_2=(0,0,1,0,\dots),\quad\dots
$$
and (cf.~\eqref{eqvFourier})
$$
\phi=(\phi_0,\phi_1,\phi_2,\dots),\qquad
v(\cdot,t)=(v_0(t),v_1(t),v_2(t),\dots).
$$

\begin{figure}[ht]
      {\ \hfill\epsfxsize70mm\epsfbox{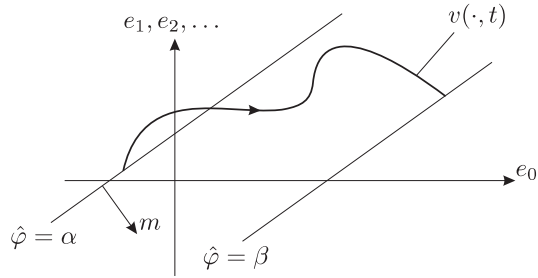}\hfill\ }
        \caption{The plane spanned by  $e_0=(1,0,0,\dots)$ and
$m=(m_0,m_1,m_2,\dots)$.}
        \label{figGeomej}
\end{figure}

 Consider the plane going through the origin and spanned
by the vector $e_0=(1,0,0,\dots)$ and the vector
$m=(m_0,m_1,m_2,\dots)$ (if they are parallel, i.e.,
$m_1=m_2=\dots=0$, then we consider an arbitrary plane containing
$e_0$). We note that the angle between the vectors $m$ and $e_0$
is acute (their scalar product is equal to $m_0>0$).
 Clearly, the orthogonal projection of the hyperspace
$\hat\phi=\sum\limits_{j=0}^\infty m_j\phi_j=\alpha$ (or $\beta$)
on this plane is a line (see Fig.~\ref{figGeomej}).

\begin{figure}[ht]
      {\ \hfill\epsfxsize45mm\epsfbox{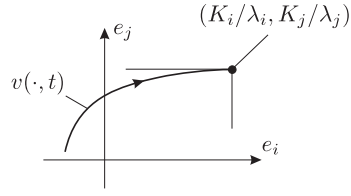}\hfill\ }
        \caption{The plane spanned by  $e_i$ and $e_j$, $i\ne j$,
        $i,j\ge 1$.}
        \label{figGeomej12}
\end{figure}

Due to~\eqref{eqExistenceUniquenessNoHyst5''}, $v_0(t)$ ``goes''
from the left to the right with the constant speed $K_0>0$, while
$v_j(t)$ exponentially converge  to $K_j/\lambda_j$ (see
Fig.~\ref{figGeomej12}).

Due to Theorem~\ref{thExistenceUniquenessNoHystMild}, the original
problem~\eqref{eq2.1}--\eqref{eq2.3} can be written as follows:
\begin{equation}\label{eqDS}
\begin{aligned}
\dot v_0(t)&=\cH(\hat v)(t)K_0, & v_0(0)&=\phi_0,\\
\dot v_j(t)&=-\lambda_j  v_j(t)+\cH(\hat v)(t)K_j, &
 v_j(0)&=\phi_j\qquad (j=1,2\dots).
 \end{aligned}
\end{equation}
Equations~\eqref{eqDS} define an infinite-dimensional dynamical
system for the functions $v_j(t)$. These functions are ``coupled''
via formula~\eqref{eqMeanTempFourier} for the mean temperature,
which is the argument of the hysteresis operator $\cH$.

\subsection{Invariant subsystem and ``guiding-guided'' decomposition}

In this subsection, we show that if some coefficients $m_j$
vanish, then the system~\eqref{eqDS} has an invariant subsystem.

We introduce the sets of indices
$$
\bbJ=\{j\in\bbN : m_j\ne 0\},\qquad \bbJ_0=\{j\in\bbN : m_j=0\}.
$$
Clearly, $\{0\}\cup\bbJ\cup\bbJ_0=\{0,1,2,\dots\}$.

Note that, for any solution $v(x,t)$ of
problem~\eqref{eq2.1}--\eqref{eq2.3}, we have
(cf.~\eqref{eqMeanTempFourier})
$$
\hat v(t)=\sum\limits_{j\in\{0\}\cup\bbJ} m_jv_j(t),\quad t\ge0.
$$
Therefore, the dynamics of $v_j(t)$, $j\in\{0\}\cup\bbJ$, does not
depend on the functions $v_j(t)$, $j\in\bbJ_0$, and is described
by the invariant dynamical system
\begin{equation}\label{eqInvariantDSN}
\begin{aligned}
\dot v_0(t)&=\cH(\hat v)(t)K_0, & v_0(0)&=\phi_0,\\
\dot v_j(t)&=-\lambda_j  v_j(t)+\cH(\hat v)(t)K_j, &
 v_j(0)&=\phi_j\qquad (j\in\bbJ).
 \end{aligned}
\end{equation}

The dynamics of $v_j(t)$, $j\in\bbJ_0$, is described by the system
\begin{equation}\label{eqInfiniteDSN}
\dot v_j(t)=-\lambda_j  v_j(t)+\cH(\hat v)(t)K_j,\qquad
 v_j(0)=\phi_j\qquad (j\in\bbJ_0),
\end{equation}
where the hysteresis operator $\cH$ depends only on the functions
$v_j(t)$ from the system~\eqref{eqInvariantDSN}.

\begin{definition}\label{defGuidingSystem}
We say that the system~\eqref{eqInvariantDSN} is {\it guiding},
while  the system~\eqref{eqInfiniteDSN} is  {\it   guided}
(by~\eqref{eqInvariantDSN}).
\end{definition}

In what follows, we will use the following notation. For any
number $\phi_0$ and (possibly, infinite-dimensional) vectors
$\{\phi_j\}_{j\in\bbJ}$ and $\{\phi_j\}_{j\in\bbJ_0}$
($\phi_j\in\bbR$), we denote
$$
\bphi =\{\phi_j\}_{j\in\bbJ},\qquad
\tilde\bphi=\{\phi_j\}_{j\in\{0\}\cup\bbJ},\qquad
\bphi_0=\{\phi_j\}_{j\in\bbJ_0}.
$$
Thus, e.g., $\tilde\bv(t)$ and $\bv_0(t)$ will represent the
solutions of the guiding system~\eqref{eqInvariantDSN} and the
guided system~\eqref{eqInfiniteDSN}, respectively.

The above decomposition of the system~\eqref{eqDS} implies the
corresponding decomposition of the phase space~$H^1$:
\begin{equation}\label{eqHDecomp}
 H^1=\bbR\times V\times V_0=\tilde V\times V_0,
\end{equation}
where the norms in $V$, $\tilde V$, and $V_0$ are given by
\begin{equation}\label{eqNormV}
\|\bphi\|_V=\left(\sum\limits_{j\in\bbJ}(1+\lambda_j)|\phi_j|^2\right)^{1/2},\quad
\|\tilde\bphi\|_{\tilde V} =
\left(\sum\limits_{j\in\{0\}\cup\bbJ}(1+\lambda_j)|\phi_j|^2\right)^{1/2},\quad
\|\bphi_0\|_{V_0}=\left(\sum\limits_{j\in\bbJ_0}(1+\lambda_j)|\phi_j|^2\right)^{1/2}.
\end{equation}

Further we  show that   Definition~\ref{defGuidingSystem} is quite
natural. In particular, we  prove that if $\tilde\bz(t)$ is a
periodic solution of the guiding system~\eqref{eqInvariantDSN},
then there exists a periodic solution of the full
system~\eqref{eqDS}  of the form $(\tilde\bz(t),\bz_0(t))$.
Moreover,  the latter is stable if and only if $\tilde\bz(t)$ is a
stable periodic solution of the guiding
system~\eqref{eqInvariantDSN}.

As an application of this result, assuming that the set $\bbJ$ is
finite, we will construct a  periodic solution $z(x,t)$ with small
period such that $\tilde\bz(t)$ is an unstable periodic solution
of the guiding system~\eqref{eqInvariantDSN}. Clearly, $z(x,t)$
will be unstable in this case, too.

%% file: sect3.tex
\section{Periodic Solutions}\label{secPeriodic}

\subsection{Conditional existence of periodic
solution}\label{subsecConditional}

We begin with a definition of periodic solutions (with two
switchings on the period) of problem~\eqref{eq2.1}, \eqref{eq2.3}.
Recall that the symbol
  $\hat{\phantom{a}}$  refers to the ``average'' of the
function (see Sec.~\ref{subsecSetting}).

\begin{definition}\label{defSTPeriodicSolution}
A function $z(x,t)$  is called an {\em    $(s,\sigma)$-periodic
solution {\rm (}with period $T=s+\sigma${\rm )} of
problem~\eqref{eq2.1}, \eqref{eq2.3}\/}  if there is a function
$\psi\in H^1$   such that the following holds:
\begin{enumerate}
\item $\hat \psi=\alpha$,
\item $z(x,t)$ is a solution  of
problem~\eqref{eq2.1}--\eqref{eq2.3}  (in $Q_\infty$) with the
initial data $\psi$,
\item there are exactly two  switching
moments   $s$ and $T$ of $\cH(\hat z)$ on the interval $(0,T]$
(such that $\hat z(s)=\beta$ and $\hat z(T)=\alpha$),
\item
 $z(x, T)=z(x,0)\ (=\psi(x))$.
\end{enumerate}
\end{definition}

\begin{definition}\label{defPeriodicTrajectory}
If $z(x,t)$  is  an      $(s,\sigma)$-periodic solution of
problem~\eqref{eq2.1}, \eqref{eq2.3} and $T=s+\sigma$, then the
sets
$$
\Gamma=\{z(\cdot,t),\ t\in[0,T]\},\qquad
\tilde\Gamma=\{\tilde\bz(t):t\in[0,T]\},\qquad
\Gamma_0=\{\bz_0(t):t\in[0,T]\}.
$$
are called the {\em trajectories} of $z(x,t)$, $\tilde\bz(t)$, and
$\bz_0(t)$, respectively.
\end{definition}

We also consider two parts of the trajectory corresponding to the
hysteresis value $\cH(\hat z)=1$ and $-1$:
\begin{equation*}\label{eqGamma12}
\Gamma_1=\{z(\cdot,t),\ t\in[0,s]\},\qquad \Gamma_2=\{z(\cdot,t),\
t\in[s,T]\}.
\end{equation*}
Similarly, one introduces the sets $\tilde\Gamma_j$ and
$\Gamma_{0j}$, $j=1,2$.

\begin{remark}
It follows from the definition of the hysteresis operator $\cH$
and from Definition~\ref{defSTPeriodicSolution} that if $z(x,t)$
is   an $(s,\sigma)$-periodic solution with period $T=s+\sigma$ of
problem~\eqref{eq2.1}, \eqref{eq2.3}, then
$$
\cH(\hat z)(t)=1,\quad t\in[0,s);\qquad \cH(\hat z)(t)=-1,\quad
t\in[s,T).
$$
\end{remark}

Further in this section, we establish the connection between
periodic solutions of the guiding system~\eqref{eqInvariantDSN}
and those of the full system~\eqref{eqDS}. The definitions of
$(s,\sigma)$-periodic solutions for the guiding
system~\eqref{eqInvariantDSN} and for the full system~\eqref{eqDS}
are analogous to Definition~\ref{defSTPeriodicSolution}.

The following theorem generalizes Theorem~4.4
in~\cite{GurJaegerSkubSIAM}, where $m_0\ne0$ and $m_1=m_2=\dots=0$
(i.e., $\bbJ=\varnothing$ and $\bbJ_0=\bbN$).

\begin{theorem}\label{tConditionalPeriodic}
Let $\tilde\bz(t)$ be an $(s,\sigma)$-periodic solution of the
guiding system~\eqref{eqInvariantDSN}. Then there exists a unique
function $\bz_0(t)$ such that $(\tilde\bz(t),\bz_0(t))$ is an
$(s,\sigma)$-periodic solution of the full system~\eqref{eqDS}
{\rm (}which generates an $(s,\sigma)$-periodic solution $z(x,t)$
of problem~\eqref{eq2.1}, \eqref{eq2.3}{\rm )}.
\end{theorem}
\proof We recall that the spaces $\tilde V$ and $V_0$ form the
decomposition of $H^1$ (cf.~\eqref{eqHDecomp}).

We introduce a nonlinear operator $\bM_T:V_0\to V_0$  as follows.
For any $\bphi_0\in V_0$, we consider the element
$(\tilde\bz(0),\bphi_0)\in H^1$. By
Theorem~\ref{tExistenceUniqueness} and the invariance of the
guiding system~\eqref{eqInvariantDSN}, there is a unique solution
of the full system~\eqref{eqDS}, which is of the form
$(\tilde\bz(t),\bv_0(t))\in H^1$. Clearly, $\bv_0(t)$ is a
solution of the guided system~\eqref{eqInfiniteDSN}. We set
$$
\bM_T(\bphi_0)=\bv_0(T),\qquad T=s+\sigma.
$$

We claim that $\bM_T$ is a contraction map. Indeed, let
$\bphi_0^1,\bphi_0^2\in V_0$ and let $\bv_0^1(t)$ and $\bv_0^2(t)$
be the corresponding solutions of the guided
system~\eqref{eqInfiniteDSN}. Since the mean temperature is
defined via $\tilde\bz(t)$ and does not depend on $\bv_0^1(t)$ and
$\bv_0^2(t)$, it follows that the difference
$\bw_0(t)=\bv_0^1(t)-\bv_0^2(t)$ satisfies the equations
$$
\dot w_j(t)=-\lambda_j w_j(t),\qquad w_j(0)=\phi_j^1-\phi_j^2
\qquad (j\in\bbJ_0).
$$
Therefore,
$$
\|\bM_T(\bphi_0^1)-\bM_T(\bphi_0^2)\|_{V_0}^2=\|\bw_0(T)\|_{V_0}^2=\sum\limits_{j\in\bbJ_0}(1+\lambda_j)e^{-2\lambda_j
T}|\phi_j^1-\phi_j^2|^2\le e^{-2\varkappa
T}\|\bphi_0^1-\bphi_0^2\|_{V_0}^2,
$$
where $\varkappa=\min\limits_{j\in\bbJ_0}\lambda_j>0$.

Thus, $\bM_T$ has a unique fixed point $\bpsi_0\in V_0$, which
yields the desired $(s,\sigma)$-periodic solution
$(\tilde\bz(t),\bz_0(t))$
 of the full system~\eqref{eqDS}.
\endproof

Further, we will  study the connection between the stability and
attractivity of solutions of the guiding system and the guided and
full systems. To do so, we need to define the Poincar\'e maps of
the respective systems.

\subsection{The Poincar\'e maps}

In this subsection, we introduce the Poincar\'e maps for the full
system~\eqref{eqDS} and for the guiding
system~\eqref{eqInfiniteDSN}. It is proved in~\cite{GurHystUncond}
that the stability of a periodic solution of the full system
follows from the stability of the corresponding fixed point of the
Poincar\'e map. Therefore, we will concentrate on the properties
of the   Poincar\'e map.

\begin{figure}[ht]
      {\ \hfill\epsfxsize140mm\epsfbox{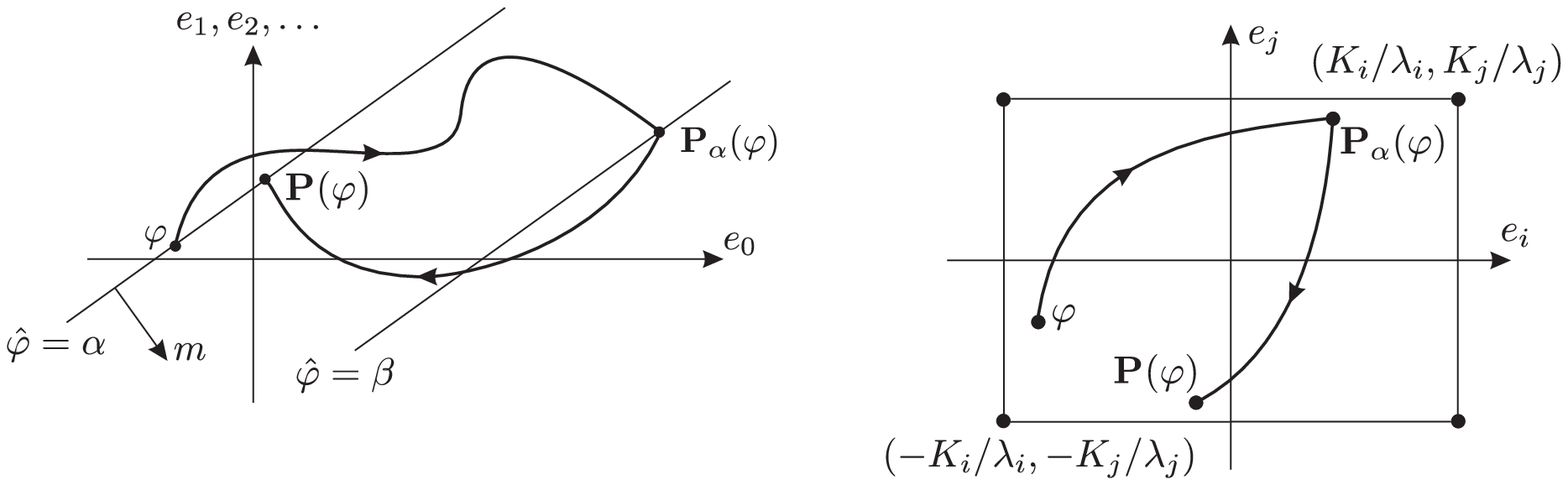}\hfill\ }
        \caption{The operators $\bP_\alpha$ and $\bP=\bP_\beta\bP_\alpha$ on the
        planes $(e_0,m)$ and $(e_i,e_j)$, $i\ne j$}
        \label{figPAlphaBeta}
\end{figure}

We consider nonlinear operators (see Fig.~\ref{figPAlphaBeta})
$$
\begin{aligned}
&\bP_\alpha:\{\phi\in H^1: \hat \phi<\beta\}\to \{\phi\in
H^1: \hat \phi=\beta\},\\
&\bP_\beta:\{\phi\in H^1: \hat \phi>\alpha\}\to \{\phi\in H^1:
\hat \phi=\alpha\}
\end{aligned}
$$
defined as follows.

Let $\phi\in H^1$, $\hat \phi<\beta$, and let $v(x,t)$ be the
corresponding solution of problem~\eqref{eq2.1}--\eqref{eq2.3} in
($Q_\infty$). Due to Theorem~\ref{tExistenceUniqueness}, there
exists the first switching moment  $t_1$   such that $\hat
v(t_1)=\beta$ and there are no other switchings on the interval
$(0,t_1)$.  In other words, the function $v^\alpha(x,t):=v(x,t)$
is a solution of the initial boundary-value problem on the
interval $(0,t_1)$:
\begin{equation}\label{eqv1_1}
v^\alpha_t(x,t) =\Delta v^\alpha(x,t) \quad ((x,t)\in Q_{t_1}),
\end{equation}
\begin{equation}\label{eqv1_2}
v^\alpha(x,0)=\phi(x)\quad (x\in Q),
\end{equation}
\begin{equation}\label{eqv1_3}
  \frac{\partial v^\alpha}{\partial\nu} =  K(x)    \quad
((x,t)\in\Gamma_{t_1}).
\end{equation}
We set $\bP_\alpha(\phi)=v^\alpha(\cdot,t_1)$.

The operator $\bP_\beta$ is defined in a similar way. Let $\phi\in
H^1$ and $\hat \phi>\alpha$. As before, there is a moment
$\tau_2>0$ and a function $v^{\beta}(x,t)$  such that
$v^{\beta}(x,t)$ is a solution of the problem
\begin{equation}\label{eqv2_1}
v_t^\beta(x,t) =\Delta v^\beta(x,t) \quad ((x,t)\in Q_{\tau_2}),
\end{equation}
\begin{equation}\label{eqv2_2}
v^\beta(x,0)=\phi(x)\quad (x\in Q),
\end{equation}
\begin{equation}\label{eqv2_3}
  \frac{\partial v^\beta}{\partial\nu} = - K(x)    \quad
((x,t)\in\Gamma_{\tau_2}),
\end{equation}
$\widehat{v^\beta}(\tau_2)>\alpha$ for $t<\tau_2$, and
$\widehat{v^\beta}(\tau_2)=\alpha$. We set
$\bP_\beta(\phi)=v^\beta(\cdot,\tau_2)$.

We introduce the {\it Poincar\'e map} for
problem~$\eqref{eq2.1}$--$\eqref{eq2.3}$, or, equivalently, for
the full system~$\eqref{eqDS}$
$$
\begin{gathered}
\bP:\{\phi\in H^1: \hat \phi<\beta\}\to \{\phi\in
H^1: \hat \phi=\alpha\},\\
\bP=\bP_\beta\bP_\alpha.
\end{gathered}
$$

We also introduce the operator (functional) $\bt_1:\{\phi\in H^1:
\hat \phi<\beta\}\to\bbR$ given by
$$
\bt_1(\phi)=\text{the first switching moment of } \cH(\hat v)\
\text{for system~\eqref{eq2.1}--\eqref{eq2.3}.}
$$

We will use the following result (see Remark 4.3
in~\cite{GurHystUncond}).

\begin{lemma}\label{lPContinuous}
Let $z(x,t)$ be a periodic solution of problem~\eqref{eq2.1},
\eqref{eq2.3}. If
$$
 \dfrac{d\hat z}{dt}\ne 0\quad\text{at the
switching moments},
$$
 then the operators
$\bP_\alpha(\phi)$, $\bP(\phi)$, and $\bt_1(\phi)$ are
continuously differentiable in a neighborhood of
$\Gamma_1\cap\{\phi\in H^1:\hat\phi<\beta\}$. The operator
$\bP_\beta(\phi)$ is continuously differentiable in a neighborhood
of $\Gamma_2\cap\{\phi\in H^1:\hat\phi>\alpha\}$.
\end{lemma}

We denote by $\tilde\bE:H^1\to\tilde V$ the orthogonal projector
from $H^1$ onto $\tilde V$.

Similarly to the operators $\bP_\alpha$, $\bP_\beta$,  $\bP$, and
$\bt_1$, we introduce the operators $\tilde\Pi_\alpha$,
$\tilde\Pi_\beta$, $\tilde\Pi$, and $\tilde\bt_1$, respectively,
corresponding to the invariant guiding
system~\eqref{eqInvariantDSN} and defined on the elements from
$\tilde V$. Due to the invariance of~\eqref{eqInvariantDSN}, we
have
$$
\begin{gathered}
\tilde\Pi_\alpha(\tilde\bphi)=\tilde\bE\bP_\alpha(\tilde\bphi,\bphi_0),\quad
\tilde\Pi_\beta(\tilde\bphi)=\tilde\bE\bP_\beta(\tilde\bphi,\bphi_0),\quad
\tilde\Pi(\tilde\bphi)=\tilde\bE\bP(\tilde\bphi,\bphi_0),\quad
\tilde\bt_1(\tilde\bphi)=\bt_1(\tilde\bphi,\bphi_0)\\
  \forall \tilde\bphi\in\tilde V,\ \forall \bphi_0\in V_0.
\end{gathered}
$$

We say that $\tilde\Pi$ is the {\it guiding Poincar\'e map}.

%
%

The following theorem shows that the stability of (exponential
convergence to) a fixed point of the guiding Poincar\'e map
$\tilde\Pi$ implies the stability of (exponential convergence to)
the fixed point of the Poincar\'e map $\bP$ of the full system.

First, we introduce some notation. Let $z(x,t)$ be an
$(s,\sigma)$-periodic solution with period $T=s+\sigma$ of
problem~\eqref{eq2.1}, \eqref{eq2.3} and $(\tilde\bz(t),\bz_0(t))$
the corresponding periodic solution of the full
system~\eqref{eqDS}. We denote
$$
\tilde\bpsi=\tilde\bz(0),\qquad \bpsi_0=\bz_0(0).
$$

Let $v(x,t)$ be another  solution of
problem~\eqref{eq2.1}--\eqref{eq2.3} such that $\hat v(0)=\alpha$,
and let $(\tilde\bv(t),\bv_0(t))$ be the corresponding  solution
of the full system~\eqref{eqDS}. The initial data will be denoted
by
$$
\tilde\bphi=\tilde\bv(0),\qquad \bphi_0=\bv_0(0)
$$
and the consecutive switching moments by $t_1,t_2,\dots$. We also
set $t_0=0$.

\begin{theorem}\label{tGuidedConverge}
Suppose that
$$
 \dfrac{d\hat z}{dt}\ne 0\quad\text{at the
switching moments}.
$$
\begin{enumerate}
\item
For any $\delta_0>0$, there exists $\tilde\delta>0$ such that if
$\tilde\bv(t_i)$ remain in the $\tilde\delta$-neighborhood of
$\tilde\bpsi$ for even $i$ and in the $\tilde\delta$-neighborhood
of $\tilde\bz(s)$ for odd~$i$ {\rm (}$i=0,1,2,\dots${\rm )}, then,
for all $\bphi_0$ in the $\delta_0$-neighborhood  of $\bpsi_0$,
$\bv_0(t_i)$  remain in the $\delta_0$-neighborhood  of $\bpsi_0$
for even $i$ and in the $\delta_0$-neighborhood of $\bz_0(s)$ for
odd $i${\rm ;}
\item Let
\begin{equation}\label{eqGuidingConverge}
\|\tilde\bv(t_i)-\tilde\bpsi\|_{\tilde
V}+\|\tilde\bv(t_{i+1})-\tilde \bz(s)\|_{\tilde V} \le \tilde
k\tilde q^i,\qquad i=0,2,4\dots,
\end{equation}
for some $0<\tilde q<1$ and $\tilde k>0$  which do not depend on
$i$. Then, for any neighborhood $\cV_0$ of $\bpsi_0$ and for all
$\bphi_0\in\cV_0$,
\begin{equation}\label{eqGuidedConverge}
\|\bv_0(t_i)-\bpsi_0\|_{V_0}+\|\bv_0(t_{i+1})-\bz(s)\|_{V_0}\le
k_0q_0^i,\qquad i=0,2,4\dots,
\end{equation}
where $0<q_0=q_0(\tilde q)<1$ and $k_0=k_0(\tilde k,\cV_0,\tilde
q)>0$ do not depend on $\tilde\bphi$ and $\bphi_0$ in the
corresponding neighborhoods.
\end{enumerate}
\end{theorem}

In the proof of this theorem, we will use the following technical
lemma.

\begin{lemma}\label{lAbstractSequence}
Let a sequence $b_0,b_2,b_4,\dots$ of nonnegative numbers
satisfies the inequalities
$$
b_{i+2}\le \zeta b_i+k\nu^i,\qquad i=0,2,4,\dots,
$$
where $k>0$ and $0<\zeta,\nu<1$ do not depend on    $i$. Then
there are numbers $0<q=q(\zeta,\nu)<1$ and $c=c(\zeta,\nu,q)>0$
which do not depend on   $i$  such that
$$
b_i\le (b_0+c)q^i,\qquad  i=0,2,4,\dots.
$$
\end{lemma}
\proof Let $\gamma=\max(\zeta,\nu^2)$. Clearly, $0<\gamma<1$. Then
$$
b_{i+2}\le \gamma b_i+k\gamma^{i/2},\qquad i=0,2,4,\dots.
$$
Consider the sequence $c_i=b_i\gamma^{-i/2}$. It  satisfies
$$
c_{i+2}\le c_i+k\gamma^{-1}\le c_0+k_1 i,
$$
where $k_1>0$ does not depend on $c_i$ and $i$, which yields the
desired estimate of $b_i$.
\endproof

\proof[Proof of Theorem~\ref{tGuidedConverge}] 1. Let us prove
assertion 1.

1a. By assumption, $s=\tilde\bt_1(\tilde\bpsi)$ is the first
switching moment of $\cH(\hat z)$. Denote by
$\tau=t_1=\tilde\bt_1(\tilde\bphi)$ the first switching moment of
$\cH(\hat v)$. By Lemma~\ref{lPContinuous},
$\tilde\bt_1(\tilde\bphi)$ is continuously differentiable in a
sufficiently small $\tilde\delta_1$-neighborhood of $\tilde\bpsi$.
Hence,
\begin{equation}\label{eqGuidedConverge1}
|\tau-s|\le k_1\|\tilde\bphi-\tilde\bpsi\|_{\tilde V}\le
k_1\tilde\delta,
\end{equation}
where $k_1>0$ depends on $\tilde\delta_1$ but does not depend on
$\tilde\delta\le\tilde\delta_1$.

1b. Using Remark~\ref{rExplicitVj} and the fact that, for any
$\varepsilon>0$, there is $d_\varepsilon>0$ such that
$$
(a+b)^2\le (1+\varepsilon)a^2+d_\varepsilon b^2,
$$
we obtain
\begin{equation}\label{eqGuidedConverge2}
\begin{aligned}
&\|\bv_0(\tau)-\bz_0(s)\|_{V_0}^2
=\sum\limits_{j\in\bbJ_0}(1+\lambda_j)
|v_j(\tau)-z_j(s)|^2\\
&\qquad=\sum\limits_{j\in\bbJ_0}(1+\lambda_j)
\left|\left(\phi_j-\dfrac{K_j}{\lambda_j}\right)e^{-\lambda_j\tau}-\left(\psi_j-
\dfrac{K_j}{\lambda_j}\right)e^{-\lambda_j s}\right|^2\\
&\qquad\le (1+\varepsilon)\sum\limits_{j\in\bbJ_0}(1+\lambda_j)
|\phi_j-\psi_j|^2 e^{-2\lambda_j \tau} +
d_\varepsilon\sum\limits_{j\in\bbJ_0}(1+\lambda_j)
\left|\psi_j-\dfrac{K_j}{\lambda_j}\right|^2
\left|e^{-\lambda_j\tau}-e^{-\lambda_j s}\right|^2.
\end{aligned}
\end{equation}

 Now we fix
$\varepsilon>0$  such  that
\begin{equation}\label{eqsigmaless1}
\zeta=(1+2\varepsilon)e^{-2\varkappa s}<1,
\end{equation}
where  $\varkappa=\min\limits_{j\in\bbJ_0}\lambda_j>0$. Further,
taking into account~\eqref{eqGuidedConverge1}, we choose
$\tilde\delta>0$ so small that
\begin{equation}\label{eqGuidedConverge2'}
(1+\varepsilon)e^{-2\lambda_j\tau}\le(1+2\varepsilon)e^{-2\varkappa
s}.
\end{equation}

Combining~\eqref{eqGuidedConverge2}, \eqref{eqsigmaless1},
and~\eqref{eqGuidedConverge2'} yields
\begin{equation}\label{eqGuidedConverge3}
\begin{aligned}
&\|\bv_0(\tau)-\bz_0(s)\|_{V_0}^2
 \le  \zeta\|\bphi_0-\bpsi_0\|_{V_0}^2 +
d_\varepsilon\sum\limits_{j\in\bbJ_0}(1+\lambda_j)
\left|\psi_j-\dfrac{K_j}{\lambda_j}\right|^2
\left|e^{-\lambda_j\tau}-e^{-\lambda_j s}\right|^2.
\end{aligned}
\end{equation}

 Using~\eqref{eqKjConvergence}, \eqref{eqNormV}, and
estimate~\eqref{eqGuidedConverge1}, we deduce
from~\eqref{eqGuidedConverge3}
$$
\|\bv_0(\tau)-\bz_0(s)\|_{V_0}^2 \le
\zeta\delta_0^2+k_2(\tilde\delta),
$$
where $k_2(\tilde\delta)>0$ and $k_2(\tilde\delta)\to 0$  as
$\tilde\delta\to 0$. In particular, this implies that
$\bv_0(t_1)=\bv_0(\tau)$ belongs to the $\delta_0$-neighborhood of
$\bz_0(s)$, provided that $\tilde\delta=\tilde\delta(\delta_0)$ is
sufficiently small.

In the same way, one can now show that $\bv_0(t_2)$ belongs to the
$\delta_0$-neighborhood of $\bpsi_0=\bz_0(T)$. By induction, we
obtain assertion 1.

2. Now we prove that
\begin{equation}\label{eqGuidedConverge4}
\|\bv_0(t_i)-\bpsi_0\|_{V_0}\le k_0q_0^i,\qquad i=0,2,4\dots.
\end{equation}
The rest part of estimate~\eqref{eqGuidedConverge} can be proved
analogously.

2a. First, we assume that $\tilde\bv(0)$ is in a sufficiently
small $\tilde\delta$-neighborhood of $\tilde\bpsi$. Then,
similarly to~\eqref{eqGuidedConverge3}, we have for even $i$
\begin{equation}\label{eqGuidedConverge5}
\begin{aligned}
&\|\bv_0(t_{i+1})-\bz_0(s)\|_{V_0}^2
 \le  \zeta\|\bv_0(t_i)-\bpsi_0\|_{V_0}^2 +
d_\varepsilon\sum\limits_{j\in\bbJ_0}(1+\lambda_j)
\left|\psi_j-\dfrac{K_j}{\lambda_j}\right|^2 \left|e^{-\lambda_j
\tau_{i+1}}-e^{-\lambda_j s}\right|^2,
\end{aligned}
\end{equation}
where $\tau_{i+1}=t_{i+1}-t_i=\tilde\bt_1(\tilde\bv(t_i))$. Using
the differentiability of $\tilde\bt_1$ and
estimate~\eqref{eqGuidingConverge}, we have
\begin{equation}\label{eqGuidedConverge6}
|\tau_{i+1}-s|\le k_1\|\tilde\bv(t_i)-\tilde\bpsi\|_{\tilde V}\le
k_3\tilde q^i
\end{equation}

Due to~\eqref{eqGuidedConverge6}, we can assume that
$\tau_{i+1}\ge s/2$. Then, taking into
account~\eqref{eqGuidedConverge6}, we have
$$
\left|e^{-\lambda_j\tau_{i+1}}-e^{-\lambda_j s}\right|\le\lambda_j
e^{-\lambda_j s/2}|\tau_{i+1}-s|\le
 k_4
\tilde q^i.
$$
Combining this inequality with~\eqref{eqGuidedConverge5},
\eqref{eqNormV}, and~\eqref{eqKjConvergence} yields
$$
\|\bv_0(t_{i+1})-\bz_0(s)\|_{V_0}^2\le
\zeta\|\bv_0(t_i)-\bpsi_0\|_{V_0}^2 + k_5\tilde q^{2i}.
$$
Making one more step and using the last inequality, we obtain
\begin{equation}\label{eqGuidedConverge7}
\begin{aligned}
\|\bv_0(t_{i+2})-\bpsi_0\|_{V_0}^2&\le \zeta
\|\bv_0(t_{i+1})-\bz_0(s)\|_{V_0}^2 + k_5\tilde q^{2i}\le
\zeta(\zeta\|\bv_0(t_i)-\bpsi_0\|_{V_0}^2 + k_5\tilde
q^{2i})+k_5\tilde q^{2i}\\
&\le \zeta^2\|\bv_0(t_i)-\bpsi_0\|_{V_0}^2+k_6 \tilde q^{2i}.
\end{aligned}
\end{equation}
Due to Lemma~\ref{lAbstractSequence}, the latter inequality
implies~\eqref{eqGuidedConverge4}.

2b. Now we take an arbitrary $\tilde\bv(0)$ in the $\tilde
k$-neighborhood of $\tilde\bpsi$. Due
to~\eqref{eqGuidingConverge}, there exists an even number~$I$
(which does not depend on $\tilde\bv(0)$) such that
$\tilde\bv(t_I)$ is in the $\tilde\delta$-neighborhood of
$\tilde\bpsi$. Then the inequality in~\eqref{eqGuidedConverge4}
holds for $i=I,I+2,I+4,\dots$ due to part 2a of the proof.

Theorem~\ref{tExistenceUniqueness} implies the existence of
$\theta>0$ (which depends on $\tilde k$ and $\cV_0$ but does not
depend on $\tilde\bv(0)$ and $\bv_0(0)$) such that $t_I\le\theta$.
Furthermore, Theorem~\ref{tExistenceUniqueness} implies that
$$
\max\limits_{t\in[0,\theta]}\|\bv_0(t)\|_{V_0}\le
k_7(\theta,\tilde k,\cV_0)= k_8(\tilde k,\cV_0).
$$
Hence, the inequality in~\eqref{eqGuidedConverge4} holds for
$i=0,2,4,\dots$.
\endproof

\begin{remark} It follows from the proof of
Lemma~\ref{lAbstractSequence} that the convergence rate $q$ is
greater than $\gamma^{1/2}=\max(\zeta^{1/2},\nu)$ but can be
chosen arbitrarily close to this number.

Therefore,  the convergence rate $q_0$ in
estimate~\eqref{eqGuidedConverge} is greater than
$\max\left(e^{-\varkappa s},e^{-\varkappa(T-s)},\tilde q \right)$
but can be chosen arbitrarily close to this number.
\end{remark}

\subsection{Conditional attraction and stability of periodic solution}

Let $z(x,t)$ and $v(x,t)$ be the same as above, but now we  do not
assume that $\hat v(0)$ is necessarily equal to $\alpha$.

The following theorem shows that the convergence to the periodic
orbit in the guiding system implies the convergence to the
corresponding periodic orbit in the full system. Thus, we call the
phenomenon in that theorem the {\it conditional attraction}.

For the trajectories, we will use the notation given in
Sec.~\ref{subsecConditional}.

\begin{theorem}\label{tCondConvergeT-T}
Suppose that
$$
 \dfrac{d\hat z}{dt}\ne 0\quad\text{at the
switching moments}.
$$
Let
$$
\begin{cases}
\dist(\tilde\bv(t),\tilde\Gamma_1)\le \tilde k \tilde q^t &
\text{if}\quad \cH(\hat v)(t)=1, \\
\dist(\tilde\bv(t),\tilde\Gamma_2)\le \tilde k \tilde q^t &
\text{if}\quad \cH(\hat v)(t)=-1,
\end{cases}
$$ for
some $0<\tilde q<1$ and $\tilde k>0$. Then, for any bounded set
$\cV_0$ in $V_0$,  there exist $0<q=q(\tilde q)<1$ and $k=k(\tilde
k,\cV_0,\tilde q)>0$ such that, for all $\bphi_0\in\cV_0$ and
$t\ge 0$,
$$
\begin{cases}
\dist(v(\cdot,t),\Gamma_1)\le kq^t & \text{if}\quad \cH(\hat
v)(t)=1, \\
\dist(v(\cdot,t),\Gamma_2)\le kq^t & \text{if}\quad \cH(\hat
v)(t)=-1,
\end{cases}
$$
\end{theorem}
\proof 1. Suppose we have shown that
\begin{equation}\label{eqCondConvergeT-T1}
\|v(\cdot,t_i)-z(\cdot,0)\|_{H^1}+\|v(\cdot,t_{i+1})-z(\cdot,s)\|_{H^1}
\le   k  q^i,\qquad i=0,2,4\dots,
\end{equation}
where $0<q=q(\tilde q)<1$ and $k=k(\tilde k,\cV_0,\tilde q)>0$.
Then, using Lemma~\ref{lPContinuous} and arguing as in the proof
of Theorem~4.3 in~\cite{GurHystUncond}, we complete the proof.

So, let us prove estimate~\eqref{eqCondConvergeT-T1}.

2. Consider the intersection of the closure of $\tilde{\Gamma}_1$
with the set $\{\tilde\bphi\in\tilde V:\hat{\phi} = \beta\}$. This
intersection consists of the single point $\tilde{\bz}(s)$, where
$s$ is the switching moment of the periodic solution.

Since $\dfrac{d \hat{z}}{ d t}\Big|_{t=s} \ne 0$, it follows from
the implicit function theorem  that there exist a number $L>0$ and
a sufficiently small number $d_0> 0$ such that if
$$
|\beta - \hat{z}(\tau)| \le d
$$
for some $d \le d_0$ and $\tau \in [0, s]$, then
$$
|\tau - s|\le Ld.
$$

3. Consider $i = 1, 3, 5, \dots$. By  assumption,  there exists
$\tau \in [0, s]$ such that
\begin{equation}\label{1.0}
\|\tilde{\bv}(t_i) - \tilde{\bz}(\tau)\|_{\tilde V} \le \tilde{k}
\tilde{q}^{t_i}.
\end{equation}
This inequality together with the Cauchy--Bunyakovskii inequality
implies that there exists a constant $C_1
> 0$ such that
\begin{equation}\label{1.1}
|\beta - \hat{z}(\tau)| = |\hat{v}(t_i) - \hat{z}(\tau)| \le C_1
\tilde{k} \tilde{q}^{t_i}.
\end{equation}

3a. First, we assume that $t_i\ge \theta$, where $\theta>0$ is so
large that $C_1 \tilde{k} \tilde{q}^{\theta}\le d_0$ ($d_0$ is the
number from part 2 of the proof). Then, due to part 2 of the
proof, we have
\begin{equation}\label{1.2}
|\tau - s|\le L C_1 \tilde{k} \tilde{q}^{t_i}.
\end{equation}

It was proved in~\cite[Lemma 4.6]{GurHystUncond} that the periodic
solution $\tilde{\bz}(t)$ is uniformly Lipschitz-continuous on
$[0,T]$, which (together with~\eqref{1.2}) implies that
\begin{equation}\label{1.3}
\|\tilde{\bz}(\tau) - \tilde{\bz}(s)\|_{\tilde V}\le C_2
\tilde{k}\tilde{q}^{t_i}
\end{equation}
for some $C_2 > 0$.

Estimates~\eqref{1.0} and~\eqref{1.3} yield
$$
\|\tilde\bv(t_i) - \tilde{\bz}(s)\|_{\tilde V} \le
\tilde{k}(1+C_2)\tilde{q}^{t_i}\quad\text{for}\ t_i\ge\theta.
$$
Taking into account~\eqref{eqt*}, we have
\begin{equation}\label{1.4}
\|\tilde\bv(t_i) - \tilde{\bz}(s)\|_{\tilde V} \le \tilde k_1
\tilde q_1^i\quad\text{for}\ t_i\ge\theta,
\end{equation}
where $\tilde k_1\ge\tilde k$ and $0<\tilde{q_1}<1$.

3b. For $t_i\le\theta$, we have
\begin{equation}\label{1.5}
\|\tilde\bv(t_i) - \tilde{\bz}(s)\|_{\tilde V} \le
\dist(\tilde\bv(t_i),\tilde\Gamma_1)\le \tilde k.
\end{equation}

Combining~\eqref{1.4} and~\eqref{1.5} yields
$$
\|\tilde\bv(t_i) - \tilde{\bz}(s)\|_{\tilde V} \le \tilde k_2
\tilde q_1^i\quad\text{for all}\ t_i.
$$

 Applying similar arguments to $\tilde{v}(t_{i+1})$ and
$\tilde{z}(0)$ and using Theorem~\ref{tGuidedConverge}, we
obtain~\eqref{eqCondConvergeT-T1}.
\endproof

Now we discuss the phenomenon of {\it conditional stability}. When
studying the stability of periodic solutions, one considers its
small neighborhood. When doing so, one has to take into account
the initial state of the hysteresis operator.

\begin{definition}\label{defStability}
An $(s,\sigma)$-periodic solution  $z(x,t)$ of
problem~\eqref{eq2.1}, \eqref{eq2.3} is   {\em  stable}
  if,
 for any neighborhoods  $\cU_1$ of $\Gamma_1$ and $\cU_2$ of $\Gamma_2$ in
$H^1$, there exist  neighborhoods $\cV_1$ of $\Gamma_1$ and
$\cV_2$ of $\Gamma_2$   in $H^1$ such that if
$$
\phi\in \cV_1,\ \hat\phi<\beta \quad \text{or}\quad \phi\in
\cV_2,\ \hat\phi\ge\beta,
$$
then the solution $v(x,t)$ of problem~\eqref{eq2.1}--\eqref{eq2.3}
in $Q_\infty$ with the initial data $\phi$ satisfies for all
$t\ge0$:
$$
\begin{cases}
v\in \cU_1 &
\text{if}\quad \cH(\hat v)(t)=1, \\
v\in \cU_2 & \text{if}\quad \cH(\hat v)(t)=-1.
\end{cases}
$$

An $(s,\sigma)$-periodic solution is {\em   unstable} if it is not
 stable.
\end{definition}

\begin{definition}\label{defExpStability}
An $(s,\sigma)$-periodic solution  $z(x,t)$ of
problem~\eqref{eq2.1}, \eqref{eq2.3} is   {\em  uniformly
exponentially  stable}
  if it is stable and
there exist  neighborhoods $\cW_1$ of $\Gamma_1$ and $\cW_2$ of
$\Gamma_2$ in $H^1$ and numbers $0<q<1$ and $k>0$ such that if
$$
\phi\in \cW_1,\ \hat\phi<\beta \quad \text{or}\quad \phi\in
\cW_2,\ \hat\phi\ge\beta,
$$
then the solution  $v(x,t)$  of
problem~\eqref{eq2.1}--\eqref{eq2.3} in $Q_\infty$ with the
initial data $\phi$ satisfies
$$
\begin{cases}
\dist(v(\cdot,t),\Gamma_1)\le kq^t & \text{if}\quad \cH(\hat
v)(t)=1, \\
\dist(v(\cdot,t),\Gamma_2)\le kq^t & \text{if}\quad \cH(\hat
v)(t)=-1
\end{cases}
$$
for all $t\ge0$ uniformly with respect to $\phi$.
\end{definition}

Let $\tilde\bz(t)$ be an $(s,\sigma)$-periodic solution of the
guiding system~\eqref{eqInvariantDSN}. Then, by
Theorem~\ref{tConditionalPeriodic}, there exists a unique function
$\bz_0(t)$ such that $(\tilde\bz(t),\bz_0(t))$ is an
$(s,\sigma)$-periodic solution of the full system~\eqref{eqDS}. We
denote by $z(x,t)$ the corresponding $(s,\sigma)$-periodic
solution of problem~\eqref{eq2.1}, \eqref{eq2.3}.

\begin{theorem}\label{tConditionalStability}
Suppose that
$$
 \dfrac{d\hat z}{dt}\ne 0\quad\text{at the
switching moments}.
$$
Then the following assertions are equivalent.
\begin{enumerate}
\item The periodic solution $z(x,t)$ of problem~\eqref{eq2.1},
\eqref{eq2.3} is stable {\rm (}uniformly exponentially stable{\rm
)}.
\item The periodic solution $\tilde\bz(t)$ of
the guiding system~\eqref{eqInvariantDSN} is stable {\rm
(}uniformly exponentially stable{\rm )}.
\item The element $\tilde\bz(0)$ is a stable {\rm (}uniformly exponentially stable{\rm
)} fixed point of the Poincar\'e map $\tilde\Pi$.
\end{enumerate}
\end{theorem}
\proof Implication $1\Rightarrow 2$ is obvious.

Implication $2\Rightarrow 3$ is proved similarly to the proof of
Theorem~\ref{tCondConvergeT-T}.

To prove implication $3\Rightarrow 1$, one should use
Lemma~\ref{lPContinuous} and Theorem~\ref{tGuidedConverge} and
argue as in the proof of Lemma~4.7 and Theorem~4.4
in~\cite{GurHystUncond}.
\endproof

%% file: sect4.tex
\section{Symmetric Periodic Solutions}\label{secPeriodicReduction}

\subsection{Preliminary considerations}

It was noted in~\cite{GurHystUncond} that any
$(s,\sigma)$-periodic solution possesses a certain symmetry,
provided that it is unique. In fact a much stronger result holds,
namely, we show that any $(s,\sigma)$-periodic solution possesses
symmetry.

We underline that the results in the previous sections did not
depend on the symmetry of periodic solutions, but the results of
this section do. In particular, by exploiting the symmetry, we
give an algorithm for finding {\it all} periodic solutions with
two switchings on the period. Using their explicit form, we will
study their stability.

\begin{definition}\label{defSTSymmetricPeriodicSolution}
An  $(s,\sigma)$-periodic solution $z(x,t)$ of
problem~\eqref{eq2.1}, \eqref{eq2.3} is called {\em symmetric} if
$z_j(0)=-z_j(s)$, $j=1,2,\dots$.
\end{definition}

\begin{lemma}\label{lOnlySymmetric}
Let $z(x,t)$ be an $(s,\sigma)$-periodic solution  of
problem~\eqref{eq2.1}, \eqref{eq2.3}. Then $s=\sigma$ and $z(x,t)$
 is symmetric.
\end{lemma}
\proof Let $\psi(x) = z(x, 0) = z(x, s+\sigma)$ and $\xi(x) = z(x,
s)$.

By Remark~\ref{rExplicitVj},
\begin{equation}\label{fly1}
\xi_0 = \psi_0 + K_0 s,
\end{equation}
\begin{equation}\label{fly2}
\xi_j = \left(\psi_j - \frac{K_j}{\lambda_j}\right)e^{-\lambda_j
s} +\frac{K_j}{\lambda_j}, \quad j \geq 1.
\end{equation}
Applying   Remark~\ref{rExplicitVj} (with $K_j$ replaced by
$-K_j$), we conclude that
\begin{equation}\label{fly3}
\psi_0 = \xi_0 - K_0 \sigma,
\end{equation}
\begin{equation}\label{fly4}
\psi_j = \left(\xi_j + \frac{K_j}{\lambda_j}\right)e^{-\lambda_j
\sigma} -\frac{K_j}{\lambda_j}, \quad j \geq 1.
\end{equation}
Equalities \eqref{fly1} and \eqref{fly3} imply that $s = \sigma$.
Summing up \eqref{fly2} and \eqref{fly4} and taking into account
that  $s = \sigma$, we see that
$$
\psi_j + \xi_j = (\psi_j + \xi_j)e^{-\lambda_j s}, \quad j \geq 1.
$$
Hence, $\xi_j = -\psi_j$, and $z(x, t)$ is symmetric.
\endproof

\begin{remark}\label{ssPeriodicSolution}
Lemma~\ref{lOnlySymmetric} shows that the period (and the second
switching time) of any $(s,\sigma)$-periodic solution is uniquely
determined by the first switching time and vice versa. Therefore,
we will say {\it ``$T$-periodic solution''} or just {\it
``periodic solution''} instead of saying ``symmetric
$(s,s)$-periodic solution with period $T=2s$''.
\end{remark}

In~\cite{GurHystUncond}, it was shown that there is a number
$\delta_1\ge0$ such that if $\beta-\alpha>\delta_1$, then there
exists a periodic solution of problem~\eqref{eq2.1},
\eqref{eq2.3}. Furthermore, there is a number
$\delta_2\ge\delta_1$ such that if $\beta-\alpha>\delta_2$, then
there exists a unique periodic solution of problem~\eqref{eq2.1},
\eqref{eq2.3}; moreover, it is stable, and is a global attractor.
Both numbers $\delta_1$ and $\delta_2$ depend  on $Q$, $m$, and
$K$.

In this section, we will   formulate a sufficient condition which
may hold for arbitrarily small $\beta-\alpha$ and still provides
the existence of (symmetric) periodic solutions. We will show that
these solutions may be both stable and unstable.


\begin{lemma}\label{lsImplies2s}
Let $z(x,t)$ be a solution  of
problem~\eqref{eq2.1}--\eqref{eq2.3} with the initial data $\psi$,
$\hat\psi=\alpha$, and let $s>0$ be the first switching moment of
$\cH(\hat z)$. If
$$
z_j(s)=-\psi_j,\quad j=1,2,\dots,
$$
then $z(x,t)$ is a $($symmetric$)$ $2s$-periodic solution of
problem~\eqref{eq2.1}, \eqref{eq2.3}.
\end{lemma}
\begin{proof}
1. First, we show that there are no switchings for $t\in(s,2s)$
and that the second switching occurs exactly for $t=2s$. To do so,
we have to show that $\hat z(t)>\alpha$, or, equivalently, $\hat
z(s)-\hat z(t)<\beta-\alpha$ for $t\in(s,2s)$. Using
Remark~\ref{rExplicitVj} (with $K_j$ replaced by $-K_j$) and the
assumption that $z_j(s)=-\psi_j$, we have for $t\in(s,2s)$
\begin{equation}\label{eqsImplies2s1}
\begin{aligned}
z_j(t)&=\left(z_j(s)+\dfrac{K_j}{\lambda_j}\right)e^{-\lambda_j
(t-s)}-\dfrac{K_j}{\lambda_j}=\left(-\psi_j+\dfrac{K_j}{\lambda_j}\right)e^{-\lambda_j
(t-s)}-\dfrac{K_j}{\lambda_j},\quad j=1,2,\dots,\\
z_0(t)&=z_0(s)-K_0(t-s)=\psi_0+2K_0s-K_0t.
\end{aligned}
\end{equation}
Therefore, taking into account~\eqref{eqMeanTempFourier}, we have
\begin{equation}\label{eqsImplies2s2}
\begin{aligned}
\hat z(s)-\hat z(t)&=m_0K_0(t-s)+\sum\limits_{j=0}^\infty
m_j\left( \psi_j-\dfrac{K_j}{\lambda_j}\right)\left(e^{-\lambda_j
(t-s)}-1\right)\\
&=m_0K_0\sigma\theta+\sum\limits_{j=0}^\infty m_j \left(
\psi_j-\dfrac{K_j}{\lambda_j}\right)\left(e^{-\lambda_j
\theta}-1\right)=\hat z(\theta)-\hat\psi,
\end{aligned}
\end{equation}
where $\theta=t-s\in(0,s)$. But  $\hat
z(\theta)-\hat\psi<\beta-\alpha$ for $\theta\in(0,s)$ and $\hat
z(s)-\hat\psi=\beta-\alpha$ (because $s$ is the first switching
moment by assumption).

2. Now we show that $z(x,2s)=\psi(x)$. Indeed,
using~\eqref{eqsImplies2s1} and the assumption that
$z_j(s)=-\psi_j$, we obtain
$$
\begin{aligned}
z_j(2s)&=-\left[\left(\psi_j-\dfrac{K_j}{\lambda_j}\right)e^{-\lambda_j
s}+\dfrac{K_j}{\lambda_j}\right]=-z_j(s)=\psi_j,\quad j=1,2,\dots,\\
z_0(2s)&=\psi_0.
\end{aligned}
$$
\end{proof}

\subsection{Construction of symmetric periodic solutions}

Lemma~\ref{lsImplies2s} allows one to explicitly find all
$(s,s)$-periodic solutions according to the following algorithm.

{\bf Step 1.} For each $s>0$, we find the (unique)
$\psi_j=\psi_j(s)$ such that $v_j(s)=-\psi_j$ for $j=1,2,\dots$,
assuming that $\cH(\hat v) \equiv 1$ on the interval $[0, s)$. To
do so, we solve the equation (cf. Remark~\ref{rExplicitVj})
$$
\left(\psi_j-\dfrac{K_j}{\lambda_j}\right)e^{-\lambda_j
s}+\dfrac{K_j}{\lambda_j}=-\psi_j,
$$
which yields
\begin{equation}\label{eqphi_js}
\psi_j=\psi_j(s)=-\dfrac{K_j}{\lambda_j}\cdot\dfrac{1-e^{-\lambda_js}}{1+e^{-\lambda_js}}.
\end{equation}

We note that $\psi_j(0)=0$ and $\psi_j(s)$ monotonically decreases
and tends to $-K_j/\lambda_j$ as $s\to+\infty$.

{\bf Step 2.} We find the (unique) $\psi_0=\psi_0(s)$ such that
$\hat\psi=\alpha$. To do so, we solve the equation
$$
m_0\psi_0+\sum\limits_{j=1}^\infty m_j\psi_j=\alpha,
$$
which yields
\begin{equation}\label{eqphi_0s}
\psi_0=\psi_0(s)=\dfrac{1}{m_0}\left(\alpha-\sum\limits_{j=1}^\infty
m_j\psi_j(s)\right).
\end{equation}

Note that the function $\psi$ with the Fourier coefficients given
by~$\eqref{eqphi_js}$ and~$\eqref{eqphi_0s}$ belongs to $H^1$.
This follows from~\eqref{eqFourierPsiH1}
and~\eqref{eqKjConvergence}.

{\bf Step 3.} If the solution\footnote{Here and further, we
sometimes write $s$ after the semicolon to explicitly indicate
that the function depends on the chosen first switching time $s$
as on a parameter.} $v(x,t)=v(x,t;s)$ of
problem~\eqref{eqExistenceUniquenessNoHyst1}--\eqref{eqExistenceUniquenessNoHyst3}
with the initial data $\psi=\psi(s)$ is such that $\cH(\hat v)$
does not switch for $t<s$ and switches at the moment $t=s$, then,
by Lemma~\ref{lsImplies2s}, there exists a $2s$-periodic solution
$z(x,t;s)$ (which coincides with $v(x,t;s)$ for $t\le s$).

The switching condition is
$$
\sum\limits_{j=0}^\infty m_jv_j(s)=\beta,
$$
or, equivalently,
\begin{equation}\label{eqFs}
F(s):=m_0K_0s+2\sum\limits_{j=1}^\infty
m_j\dfrac{K_j}{\lambda_j}\cdot\dfrac{1-e^{-\lambda_js}}{1+e^{-\lambda_js}}=\beta-\alpha.
\end{equation}

To check that the switching does not occur before $s$, we note
that, due to Remark~\ref{rExplicitVj}, the mean temperature $\hat
v(t;s)$ corresponding to the initial condition~\eqref{eqphi_js},
$\eqref{eqphi_0s}$ is given by
$$
\hat v(t;s)=\sum\limits_{j=0}^\infty
m_jv_j(t;s)=\alpha+m_0k_0t+2\sum\limits_{j=1}^\infty
m_j\dfrac{K_j}{\lambda_j}\cdot\dfrac{1-e^{-\lambda_jt}}{1+e^{-\lambda_js}}.
$$
Therefore, the condition $\hat v(t;s)=\beta$ is equivalent to
$$
m_0k_0t+2\sum\limits_{j=1}^\infty
m_j\dfrac{K_j}{\lambda_j}\cdot\dfrac{1-e^{-\lambda_jt}}{1+e^{-\lambda_js}}=\beta-\alpha.
$$
Taking into account  equality~\eqref{eqFs}, we see that the
condition $\hat v(t;s)=\beta$ is equivalent to the following:
\begin{equation}\label{eqHts} H(t,s):=m_0k_0(t-s)+2\sum\limits_{j=1}^\infty
m_j\dfrac{K_j}{\lambda_j}\cdot\dfrac{e^{-\lambda_js}-e^{-\lambda_jt}}{1+e^{-\lambda_js}}=0.
\end{equation}
Moreover, the fulfillment of the inequality $H(t,s)<0$ for all
$t\in (0,s)$ is   necessary and sufficient   for the absence of
switching moments before the time moment $s$.

\begin{definition}
We will say that $F(s)$ and $H(t,s)$ are the {\it first} and the
{\it second characteristic functions}, while~\eqref{eqFs}
and~\eqref{eqHts} are the {\it first} and the {\it second
characteristic equations}, respectively.
\end{definition}

The first and the second characteristic equations will play a
fundamental role in the description of periodic solutions and
their bifurcation sets (see Theorems~\ref{tExistSymmetricSolution}
and~\ref{tAppearancePer} below).

The following lemmas describe some properties of the
characteristic functions.

\begin{lemma}\label{lFAnal}
\begin{enumerate}
\item $F(s)$ is continuous for  $s\ge 0$ and analytic for $s>0$,

\item $F(0)=0$, $F(s)$ increases for all sufficiently large $s>0$,
and $\lim\limits_{s\to+\infty}F(s)=+\infty$,

\item for each $\beta-\alpha>0$, the first characteristic
equation~$\eqref{eqFs}$ has  finitely many roots,

\item the positive zeroes of $F(s)$ are isolated and  may
accumulate only at the origin.
\end{enumerate}
\end{lemma}
\proof 1. The series in~\eqref{eqFs} is absolutely and uniformly
convergent for $\Re s\ge0$  due to the Cauchy--Bunyakovskii
inequality and~\eqref{eqKjConvergence}. Therefore, $F(s)$ is
continuous for  $s\ge 0$ and analytic for $s>0$.

Assertion  2 is now straightforward.

To prove assertion 3, we note that, for $\beta-\alpha>0$, the
(positive) roots of the first characteristic
equation~$\eqref{eqFs}$ cannot accumulate at the origin. This
follows by the continuity and the relation $F(0)=0$. The roots
cannot accumulate at infinity either (due to the monotonicity for
large $s$). Therefore, all the roots belong to a compact separated
from the origin. Now the analyticity for $s>0$ implies assertion
3.

Assertion 4 follows from the analyticity of $F(s)$ for $s>0$ and
from the monotonicity for large $s$.
\endproof

Similarly, one can prove the following lemma.
\begin{lemma}\label{lHCont}
\begin{enumerate}
\item $H(t,s)$ is continuous for  $s\ge 0$, $0\le t\le s$,

\item for each $s>0$, $H(t,s)$ is analytic in $t$ for $t>0$,

\item $H(0,s)=-F(s)$ and $H(s,s)\equiv 0$,

\item if $s>0$ and $F(s)>0$, then the second characteristic
equation~\eqref{eqHts} has no more than finitely many roots in $t$
for $t\in(0,s)$.
\end{enumerate}
\end{lemma}

Taking into account Lemmas~\ref{lOnlySymmetric}, \ref{lFAnal},
and~\ref{lHCont}, we formulate the above  algorithm as the
following theorem (also mind Remark~\ref{ssPeriodicSolution}).

\begin{theorem}\label{tExistSymmetricSolution}
\begin{enumerate}
\item For a given $\beta-\alpha>0$, there are no more than
finitely many periodic solutions of problem~\eqref{eq2.1},
\eqref{eq2.3}, which we denote $z^{(1)},\dots,z^{(N)}$. \item All
the periodic solutions $z^{(1)},\dots,z^{(N)}$  are symmetric.
\item If $s_1,\dots,s_N$ are half-periods of
$z^{(1)},\dots,z^{(N)}$, respectively, then $s_1,\dots,s_N$ are
the roots of the first characteristic equation~\eqref{eqFs}. \item
Let $s_{N+1},\dots,s_{N_1}$ be positive roots of the first
characteristic equation~\eqref{eqFs} different from
$s_1,\dots,s_N$. Then
\begin{enumerate}
\item $H(t,s_j)<0$ for all $t\in(0,s_j)$ if $j=1,\dots,N$,

\item $H(t;s_j)=0$ for some $t\in(0,s_j)$ if $j=N+1,\dots,N_1$.
\end{enumerate}
\end{enumerate}
\end{theorem}

In particular, Theorem~\ref{tExistSymmetricSolution} implies that
a positive root $s_j$ of the first characteristic
equation~\eqref{eqFs} ``generates'' a $2s_j$-periodic solution if
and only if $H(t,s_j)<0$ for all $t\in(0,s_j)$.

Now we will keep the domain $Q$ and the functions $m(x)$ and
$K(x)$ fixed, while allow the thresholds $\alpha$ and $\beta$
vary. We will classify the existence of all periodic (i.e.,
$(s,s)$-periodic) solutions with respect to the parameter $s$ and
with respect to the parameter $ \beta-\alpha$. By the existence of
a periodic solution for a given $s>0$ we mean that there exist
numbers $\alpha<\beta$ (depending on $s$) such that
problem~\eqref{eq2.1}, \eqref{eq2.3} with these $\alpha$ and
$\beta$ admits an $(s,s)$- or, equivalently, a $2s$-periodic
solution.

First, we show that one can divide the positive $s$-semiaxis into
intervals (whose union is denoted by~$L$) in the following way.
For every interval $L'\subset L$, either there are no
$2s$-periodic solutions for all $s\in L'$ or there is exactly one
$2s$-periodic solution for every $s\in L'$, which smoothly depends
on $s$ in $L'$. The complement $S$ of the union $L$ of all those
intervals will consist of points of possible bifurcation with
respect to $s$ (half-period). It will be a compact set. Typically,
$S$ will consist of finitely many points (see
Examples~\ref{exClassification}).

The compact set $\Sigma=F(S)$ will consist of points of possible
bifurcation with respect to the parameter $\beta-\alpha$. This set
divides the positive $(\beta-\alpha)$-semiaxis into open intervals
(whose union is denoted by $\Lambda$). For $\beta-\alpha$ in an
interval $\Lambda'\subset\Lambda$, the number of periodic
solutions remains constant and they smoothly depend on
$\beta-\alpha\in \Lambda'$ (see Example~\ref{exClassification}).

First, we introduce the set
$$
S_0=\{s>0: F(s)=0\}.
$$
Due to Lemma~\ref{lFAnal}, the set $S_0$ consists of no more than
countably many points, which may accumulate only at the origin.

To introduce the next set, we denote for $s>0$
\begin{equation}\label{eqtaus}
\tau(s)=\{t\in(0,s): H(t,s)=0\}.
\end{equation}
By Lemma~\ref{lHCont}, $\tau(s)$ consists of finitely many roots
of the equation $H(t,s)=0$ on the interval $t\in(0,s)$, provided
that $F(s)>0$.

Consider the set
$$
S_1=\{s>0: F(s)> 0,\ \tau(s)=\varnothing\ \text{and}\
H_t(t,s)|_{t=s}=0\}.
$$
Thus, $S_1$ consists of those $s$ for which the corresponding
trajectory $v(x,t;s)$ intersects the hyperplane $\hat\phi=\beta$
for the first time at the moment $s$ and touches it
nontransversally at this moment. Note that any number  $s\in S_1$
generates a $2s$-periodic solution.

 Consider the set
$$
S_2=\{s>0: F(s)> 0, \tau(s)\ne\varnothing,\ \text{and}\
H_t(t,s)|_{t=t'}=0\ \forall t'\in\tau(s)\}.
$$
Thus, $S_2$ consists of those $s$ for which the corresponding
trajectory $v(x,t;s)$ intersects the hyperplane $\hat\phi=\beta$
for the first time before the moment $s$ and touches it
nontransversally at each of the intersection moments (before $s$).
None of the numbers $s\in S_2$   generate a $2s$-periodic
solution.

We also introduce the set
$$
S_3=\{s>0: F(s)> 0\ \text{and}\ F'(s)=0\}.
$$
We note that the set $S_3$ consists of no more than countably many
isolated points which may accumulate only at the origin. This
follows from the analyticity of $F(s)$ for $s>0$ and from the
monotonicity for large $s$.

Now we set
$$
\qquad L=(0,\infty)\setminus \overline{S_0\cup S_1\cup S_2}.
$$
and
$$
\Sigma=F({S_1\cup S_2\cup S_3}),\qquad
\Lambda=(0,\infty)\setminus\overline{\Sigma}.
$$
We note that the above sets $S_i,L$ and $\Sigma,\Lambda$ do not
depend on $s$ or $\beta-\alpha$. They only depend on $m_j$, $K_j$,
and $\lambda_j$. We also note that the sets $S_0,\dots S_3$ and
$\Sigma$ are bounded. Indeed, $S_0$ and $S_3$ are bounded because
$F(s)$ monotonically increases for sufficiently large $s$.
Furthermore, it is proved in~\cite{GurHystUncond} that, for
sufficiently large $\beta-\alpha$ (hence for sufficiently large
$s$), the first switching moment for $v(x,t;s)$ is equal to $s$
and $\dfrac{d \hat v(t;s)}{dt}\Big|_{t=s}>0$. Therefore, $S_2$ and
$S_3$ are also bounded. The boundedness of $S_0,\dots,S_3$ implies
the boundedness of $\Sigma$.

\begin{theorem}\label{tAppearancePer}
\begin{enumerate}
\item Let $L'$ be an open interval in $L$. Then either there are
no $2s$-periodic solutions for all $s\in L'$ or, for any $s\in
L'$, there is a unique $2s$-periodic solution $z(x,t;s)$ of
problem~\eqref{eq2.1}, \eqref{eq2.3}. Moreover, the initial value
$z(x,0;s)$ smoothly depends on $s\in L'$ $($in the
$H^1$-topology$)$.

\item Let $\Lambda'$ be an open interval in $\Lambda$. Then the
number of periodic solutions of problem~\eqref{eq2.1},
\eqref{eq2.3} remains constant for all $\beta-\alpha\in\Lambda'$.
The initial values of those solutions and the first switching
times continuously depend on $\beta-\alpha\in\Lambda'$ $($in the
$H^1$-topology$)$.
\end{enumerate}
\end{theorem}
\proof 1. Let $L'$ be an open interval in $L$. For any $s\in L'$,
we denote by $\hat v(t;s)$  the mean temperature corresponding to
the initial condition~\eqref{eqphi_js}, $\eqref{eqphi_0s}$. We
recall that
$$
\hat v(t;s)=\beta
$$
if and only if
$$
H(t,s)=0.
$$

 Fix an arbitrary $s'\in L'$. Then $s'\notin S_0$, i.e.,
$F(s')\ne0$.  If $F(s')<0$, then $F(s)<0$ for all $s\in L'$
(otherwise, $F(s)=0$ for some $s\in L'$, but then $s\in S_0$,
which is impossible). In this case, every $s\in L'$ does not
generate a periodic solution.

Assume that $F(s)>0$.

 Consider the sets $\tau(s)$ given by~\eqref{eqtaus} for
$s\in L'$. We claim that if $\tau(s')=\varnothing$, then
$\tau(s)=\varnothing$ in a sufficiently small neighborhood of
$s'$; if $\tau(s')\ne\varnothing$, then $\tau(s)\ne\varnothing$ in
a sufficiently small neighborhood of $s'$. Indeed:

\begin{enumerate}
\item[1a.] Let $\tau(s')=\varnothing$. Suppose that there is a
sequence $s_i$ converging to $s$ and a sequence $t_i\in(0,s_i)$
such that $H(t_i,s_i)=0$. Taking a subsequence if needed, we can
assume that $t_i\to t'\in(0,s']$. Thus, by continuity of $H(t,s)$,
we have
\begin{equation}\label{eqAppearancePer1}
H(t',s')=0
\end{equation}

 Since $\tau(s')=\varnothing$ and $s'\notin S_1$, we have $H_t(t,s')|_{t=s'}\ne 0$.
 Therefore, by the implicit function theorem and by the identity $H(s,s)\equiv
0$, it follows that, in a neighborhood of the point $(s',s')$, the
only root (in $t$) of the equation $H(t,s)=0$ is $t=s$. Hence, all
$t_i$ lie outside a fixed neighborhood of $s'$, which means that
$t'<s'$. Together with~\eqref{eqAppearancePer1}, this yields
$\tau(s')\ne\varnothing$. This contradiction proves that
$\tau(s)=\varnothing$ in a sufficiently small neighborhood of
$s'$.

\item[1b.] Now let $\tau(s')\ne\varnothing$. Since $s'\notin S_2$,
there is $t'<s'$ such that $H(t',s')=0$ and  $H_t(t,s')|_{t=t'}\ne
0$. By  the implicit function theorem the equation $H(t,s)=0$
admits a solution $t=t(s)$  in a neighborhood of $s'$ such that
$t'=t(s')$. By regularity, $t(s)<s$ if the neighborhood is small
enough. Therefore, $\tau(s)\ne\varnothing$ in a sufficiently small
neighborhood of $s'$.
\end{enumerate}

To complete the proof of assertion 1, we choose an arbitrary
compact interval in $L'$, cover each point of it by the above
neighborhood and take a finite subcovering.

The smooth dependence of the initial value of the periodic
solution on $s\in L'$ follows from the explicit
formulas~\eqref{eqphi_js} and $\eqref{eqphi_0s}$.

2. Let $\Lambda'$ be an open interval in $\Lambda$.

Fix an arbitrary $b'\in \Lambda'$. Since $b>0$, Lemma~\ref{lFAnal}
implies that the first characteristic equation $F(s)=b'$ has
finitely many (say, $N_1$) positive roots $s_1',\dots,s_{N_1}'$.
Since $b'\notin \Sigma$, it follows that $s_j'\notin S_3$, i.e.,
$F'(s_j')\ne 0$. Therefore, for $b$ in a neighborhood of $b'$,
there exist exactly $N_1$ positive roots
$s_1=s_1(b),\dots,s_{N_1}=s_{N_1}(b)$ of the first characteristic
equation $F(s)=b$, which smoothly depend on $b$.

Further, we assume that there are $N$ ($N\le N_1$) numbers
$s_1',\dots,s_N'$ for which the minimal root of the equation
$H(t,s_j')=0$ on the interval $(0,s_j')$ is equal to $s_j'$. As
before, this means that $s_j'$ generate $2s_j'$-periodic solutions
for $j=1,\dots,N$ and do not generate periodic solutions for
$j=N+1,\dots,N_1$ (cf. Theorem~\ref{tExistSymmetricSolution}).

Since $b'>0$ and $b'\notin \Sigma$, it follows that $s_j'\notin
S_0\cup S_1\cup S_2$ ($j=1,\dots,N_1$). Therefore, similarly to
part~1 of the proof, for all $b$ in a neighborhood of $b'$, the
numbers $s_j=s_j(b)$ generate $2s_j$-periodic solutions for
$j=1,\dots,N$ and do not generate periodic solutions for
$j=N+1,\dots,N_1$.

To complete the proof of assertion 2, we choose an arbitrary
compact interval in $\Lambda'$, cover each point $b'$ of it by the
above neighborhood and take a finite subcovering.
\endproof

\begin{remark}
 Theorem~\ref{tAppearancePer} indicates the ways a new
periodic solution may appear or an existing periodic solution may
disappear, i.e., bifurcation occurs.

When varying the parameter $s$,  bifurcation may occur only if
$s\in S_0\cup S_1\cup S_2$.
\begin{enumerate}
\item   The condition $s\in S_0$ implies that $\alpha$ and $\beta$
coalesce.

\item The condition $s\in S_1$  corresponds to the tangential
approach of the   trajectory $v(x,t;s)$ to the hyperplane
$\hat\phi=\beta=\alpha+F(s)$. At the point $s$, the periodic
solution exists.  In the literature on switching (or hybrid)
systems, such a bifurcation is usually called ``grazing
bifurcation''. The corresponding Poincar\'e map will be
discontinuous at this point.

\item  The condition $s\in S_2$ also corresponds to the tangential
approach of the   trajectory $v(x,t;s)$ to the hyperplane
$\hat\phi=\beta=\alpha+F(s)$. However, at the point $s$, the
periodic solution does not exists. The switching occurs before the
trajectory comes in the ``symmetric'' position. This bifurcation
can also be called ``grazing bifurcation''.

\end{enumerate}

 When varying the parameter $\beta-\alpha>0$, bifurcation may
occur if a point $s\in F^{-1}(\beta-\alpha)$  belongs to $S_1$,
$S_2$, or $S_3$.

Grazing bifurcation occurs on $S_1$ and $S_2$ as described above.

If $s\in S_3\setminus (S_1\cup S_2)$, then  a new root of the
first characteristic equation~\eqref{eqFs} may appear and then
split into two roots (or   two existing roots may merge into one
and then disappear) as  $\beta-\alpha$ crosses the value $F(s)$.
If the first switching moment for $v(x,t;s)$ is equal to $s$
(i.e., $H(t,s)<0$ for $t<s$ or, equivalently,
$\tau(s)=\varnothing$), then a new periodic solution will appear
and then split into  two (or the two existing periodic solutions
will merge into one and then disappear). This corresponds to a
fold bifurcation.

On the other hand, if   the first switching moment for $v(x,t;s)$
is less than $s$ (i.e., $H(t,s)=0$ for some $t<s$ or,
equivalently, $\tau(s)\ne\varnothing$), then no bifurcation
happens.
\end{remark}

We consider an example  illustrating
Theorems~\ref{tExistSymmetricSolution} and~\ref{tAppearancePer}.

\begin{example}\label{exClassification} Let $Q$ be a one-dimensional domain, e.g., $Q=(0,\pi)$,
cf.~\cite{GlSpr-SIAM,GlSpr-JIntEqu,FriedmanJiang-CPDE,Pruess,GoetzHoffmannMeirmanov}.
Let the boundary condition~\eqref{eq2.3} be given by
$$
v_x(0,t)=0,\qquad v_x(\pi,t)=\cH(\hat v)(t).
$$
  From the physical point of view,
these boundary conditions model a thermocontrol process in a rod
with heat-insulation on one end and a heating (cooling) element on
the other.

It is easy to find that
$$
\lambda_0=0,\quad e_0=\sqrt{\dfrac{1}{\pi}},\quad K_0= e_0(\pi)=
\sqrt{\dfrac{1}{\pi}},
$$
$$
\lambda_j=j^2,\quad e_j(x)=\sqrt{\dfrac{2}{\pi}}\cos{jx},\quad
K_j= e_j(\pi)=(-1)^j \sqrt{\dfrac{2}{\pi}},\quad j=1,2,\dots.
$$

Let   $m_0=2$, $m_1=m_2=4$, and $m_3=m_4=\dots=0$. Then the
bifurcation diagram is depicted in Fig.~\ref{figBifm0=2}.
\begin{figure}[ht]
      {\ \hfill\epsfxsize150mm\epsfbox{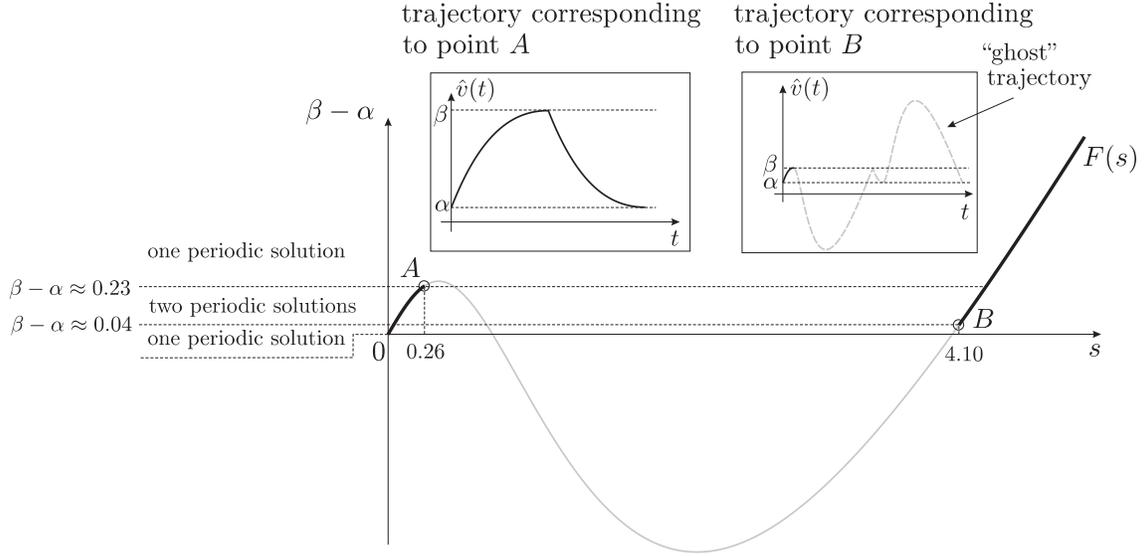}\hfill\ }
        \caption{Bifurcation diagram for $m_0=2$, $m_1=m_2=4$, and $m_3=m_4=\dots=0$. For any $s>0$, there
        exists a unique $2s$-periodic solution if the graph of $F$ is bold at the point $s$ and there are no $2s$-periodic
        solutions otherwise. For any $\beta-\alpha>0$, there exist one or two  periodic solutions depending on whether
         the horizontal
        line  levelled at $\beta-\alpha$ intersects the bold part of the graph of $F$ at one or two points,
        respectively.
        Point $A$ ($s\approx 0.26$, $\beta-\alpha\approx 0.23$) on the graph corresponds to $s\in S_1$. There exists a corresponding $2s$-periodic solution, whose
        trajectory is tangent to the hyperplane $\hat\phi=\beta$ at the moment $s$ (see the left inset). Point $B$ ($s\approx 4.10$, $\beta-\alpha\approx 0.04$)
        on the graph
        corresponds to $s\in S_2$; there does not exist a $2s$-periodic solution for this $s$. However, if one did not switch when $\hat v(t;s)$
        tangentially intersected the hyperplane $\hat\phi=\beta$ at the moment $s$, but switched only when $\hat v(t;s)$
        intersected
        the hyperplane $\hat\phi=\beta$ for the second time (at some moment $s_1>s$), then the resulting trajectory would be $2s_1$ periodic. Such a
        trajectory is referred to as a ``ghost'' trajectory (see the right
        inset).}
        \label{figBifm0=2}
\end{figure}

Let   $m_0=3.2$, $m_1=m_2=4$, and $m_3=m_4=\dots=0$. Then the
bifurcation diagram is depicted in Fig.~\ref{figBifm0=3-2}.
\begin{figure}[ht]
{\ \hfill\epsfxsize150mm\epsfbox{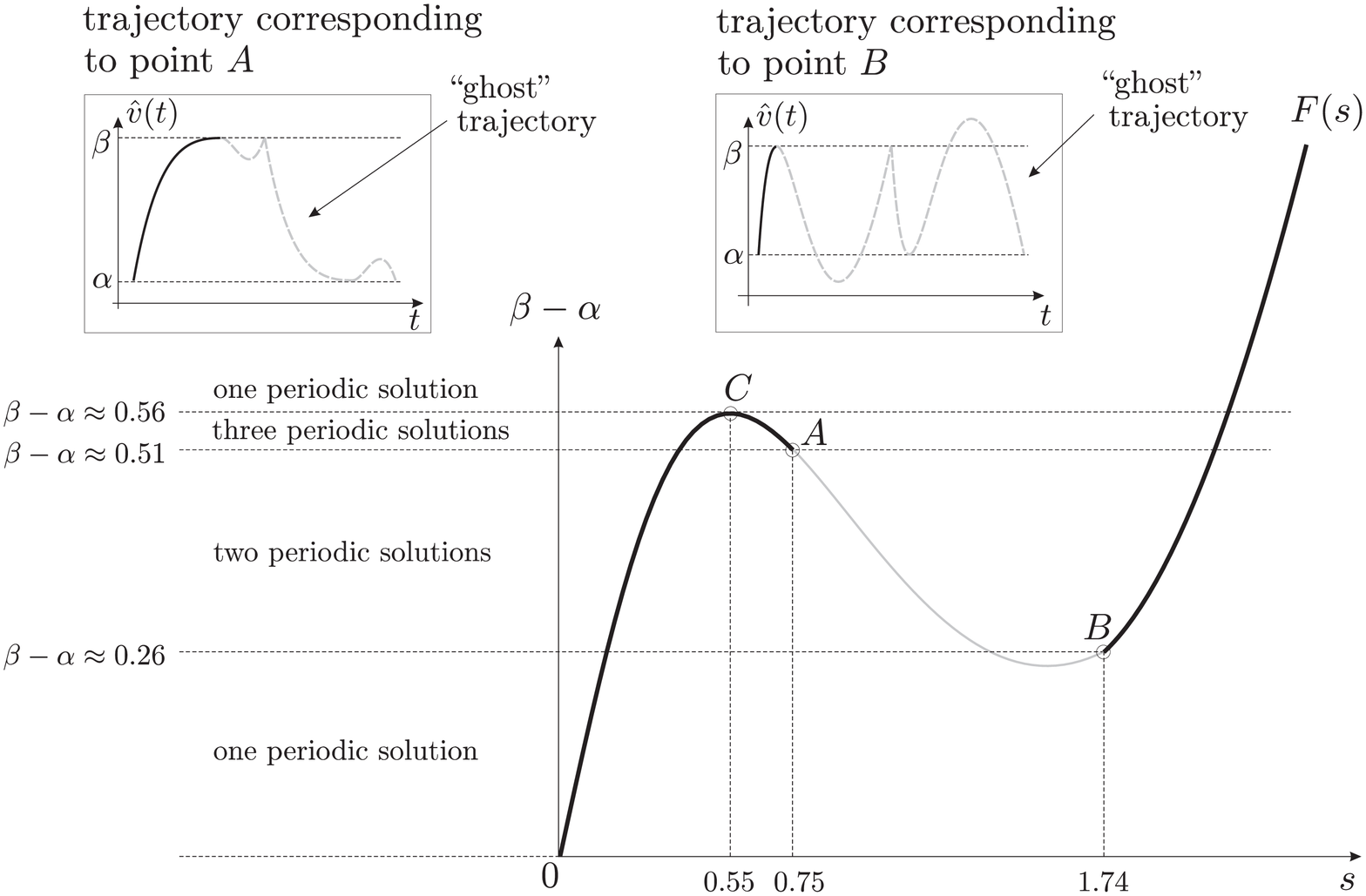}\hfill\ }
        \caption{Bifurcation diagram for $m_0=3.2$, $m_1=m_2=4$, and $m_3=m_4=\dots=0$. For any $s>0$, there
        exists a unique $2s$-periodic solution if the graph of $F$ is bold at the point $s$ and there are no $2s$-periodic
        solutions otherwise. For any $\beta-\alpha>0$, there exist one, two, or three periodic solutions depending on whether
         the horizontal
        line  levelled at $\beta-\alpha$ intersects the bold part of the graph of $F$ at one, two, or three points,
        respectively.
        Points $A$ ($s\approx 0.75$, $\beta-\alpha\approx 0.51$) and $B$ ($s\approx 1.74$, $\beta-\alpha\approx 0.26$)
        on the graph correspond to $s\in S_2$. In each of these points, a  $2s$-periodic solution does not exist.
        However, if one did not switch when $\hat v(t;s)$
        tangentially intersected the hyperplane $\hat\phi=\beta$ at the moment $s$, but switched only when $\hat v(t;s)$
        intersected
        the hyperplane $\hat\phi=\beta$ for the second time (at some moment $s_1>s$), then the resulting trajectory would be $2s_1$ periodic.
         Such a
        trajectory is referred to as a ``ghost'' trajectory (see the
        insets). Point $C$ ($s\approx 0.55$, $\beta-\alpha\approx 0.56$) corresponds to the fold bifurcation, where two periodic solutions merge
         into one and disappear as
        $\beta-\alpha$ increases and crosses the critical value $\approx 0.56$.}
        \label{figBifm0=3-2}
\end{figure}
 ``Evolution'' of periodic solutions with respect to the
parameter $\beta-\alpha$ is visualized in Fig.~\ref{figPilotka}.
\begin{figure}[ht]
      {\ \hfill\epsfxsize40mm\epsfbox{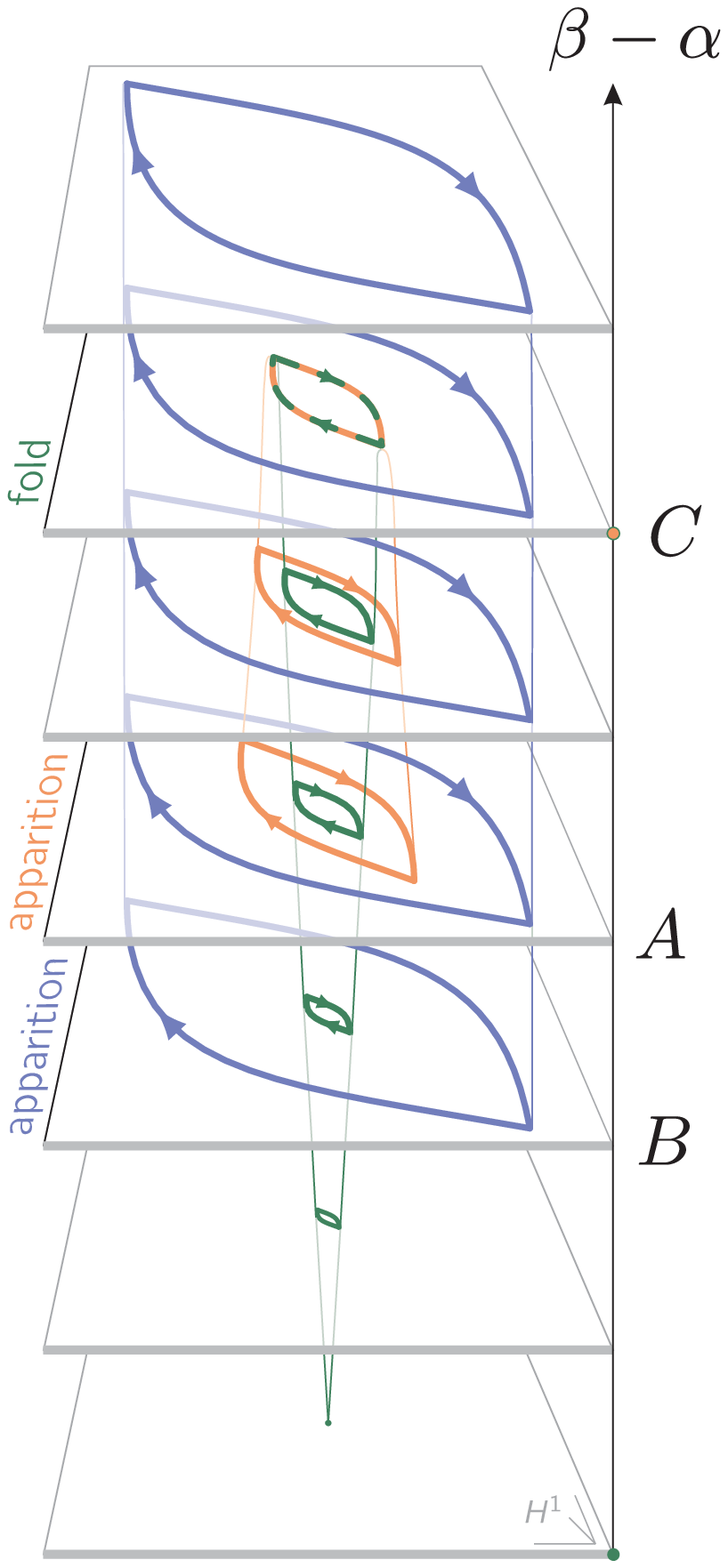}\hfill\ }
        \caption{Visualization of ``evolution'' of
periodic solutions with respect to the parameter $\beta-\alpha$
for $m_0=3.2$, $m_1=m_2=4$, and $m_3=m_4=\dots=0$. For each
$\beta-\alpha$, the horizontal plane represents the phase space
$H^1$ with periodic solutions.  Points $A$ and $B$ correspond to
apparition (or termination) of periodic solutions, while point $C$
corresponds to the fold bifurcation (cf.
Fig.~\ref{figBifm0=3-2}).}
        \label{figPilotka}
\end{figure}

%
%
\end{example}

In~\cite{GurHystUncond}, it was shown that there exists a unique
 periodic solution if $\beta-\alpha$ is large enough.
Moreover, it is stable and is a global attractor. To conclude this
section, we  prove that a periodic solution can also exist for
arbitrarily small $\beta-\alpha$. Further, we will show that such
a solution need not be stable.

Assume that the following condition holds.

\begin{condition}\label{condmKPositive}
The functions $m\in H^1$ and $K\in H^{1/2}$ satisfy
$$
M:=\sum\limits_{j=0}^\infty m_j K_j=\int_{\partial Q}
m(x)K(x)\,d\Gamma>0.
$$
\end{condition}
The convergence of the sum follows from Remark~\ref{remmjKj}. The
equality follows from the definition of $m_j$ and~$K_j$. The
essential requirement of Condition~\ref{condmKPositive} is the
positivity of the sum, or, equivalently, of the integral. From the
physical viewpoint, this condition implies the presence of thermal
sensors on a part of the boundary where the heating elements are.

\begin{theorem}\label{tExistPeriodicSolution}
Let Condition~$\ref{condmKPositive}$ hold. Then there exist
numbers $\omega>0$ and $\sigma>0$ such that, for any
$\beta-\alpha\le\omega$, there exists a  $2s$-periodic solution
$z(x,t)=z(x,t;s)$ of problem~\eqref{eq2.1}, \eqref{eq2.3} such
that $s\le\sigma$. On the interval $(0,\omega]$, the function
$s=s(\beta-\alpha)$ is strictly monotonically increasing and
$s\to0$ as $\beta-\alpha\to 0$.
\end{theorem}
\proof 1. By Condition~\ref{condmKPositive},
$F'(0)=\sum\limits_{j=0}^\infty m_j K_j>0$. Therefore, for
sufficiently small $\beta-\alpha>0$, the equation
$F(s)=\beta-\alpha$ has a unique solution $s>0$ in a small
right-hand side neighborhood $(0,\sigma]$ of the origin. Clearly,
the function $s=s(\beta-\alpha)$ possesses the properties from the
theorem.

Consider the solution $v(x,t)=v(x,t;s)$ of
problem~\eqref{eq2.1}--\eqref{eq2.3} with the initial data
$\psi=\psi(s)$ defined in Steps~1--3 above.

To complete the proof, it remains to show that $\hat v(t)=\hat
v(t;s)<\beta$ for $t<s$ and apply
Theorem~\ref{tExistSymmetricSolution}.

2. Using representation~\eqref{eqMeanTempFourier},
Remark~\ref{rExplicitVj}, and formulas~\eqref{eqphi_js}, we have
for $t\le s$
\begin{equation}\label{eqdvts}
\dfrac{d\hat v(t;s)}{dt} =m_0K_0+\sum\limits_{j=1}^\infty
m_jK_j\dfrac{2e^{-\lambda_j t}}{1+e^{-\lambda_j s}}
=M+\sum\limits_{j=1}^\infty m_jK_j\left(\dfrac{2e^{-\lambda_j
t}}{1+e^{-\lambda_j s}}-1\right).
\end{equation}
Using Remark~\ref{remmjKj}, one can easily check that the absolute
value of the series on the right-hand side is less than $M/2$ for
sufficiently small $s$ and $t\le s$. Therefore,  $\hat v(t;s)$ is
monotonically increasing until the first switching moment. Thus,
the first switching occurs for $t=s$.
\endproof

We stress that Theorem~\ref{tExistPeriodicSolution} ensures the
uniqueness of a  periodic solution with a small first switching
time~$s$ (hence small $\beta-\alpha$). However, the theorem does
not forbid the existence of other  periodic solutions with large
period and large $\beta-\alpha$.

\begin{remark} It is an open question whether one can choose the
 functions $m(x)$ and $K(x)$ and the parameters $\alpha$ and $\beta$ in
such a way that problem~\eqref{eq2.1}, \eqref{eq2.3} has no
periodic solutions.
\end{remark}

\subsection{Stability of  periodic solutions}
In this section, we will show that the thermocontrol problem with
hysteresis may admit unstable  periodic solutions.

For simplicity, we assume that only finitely many Fourier
coefficients $m_j$ do not vanish (but see
Remark~\ref{remNinfinity}).
\begin{condition}\label{condmN0}
There is $N\ge1$ such that
$$
\bbJ=\{m_0,m_1,\dots,m_N\}.
$$
\end{condition}
Clearly, modifications needed if $\bbJ$ consists of other Fourier
coefficients $m_j$ are trivial.

\begin{remark}
The fulfilment of Condition~\ref{condmN0} implies that $m\in H^1$.
Moreover, the sum in Condition~\ref{condmKPositive} becomes
finite:
$$
\sum\limits_{j=0}^N m_jK_j>0.
$$
\end{remark}

\begin{remark}\label{remN=0}
If $N=0$, i.e., $\bbJ=\{m_0\}$, then it is easy to see that the
(one-dimensional) guiding system~\eqref{eqInvariantDSN} has a
unique periodic solution for any $\alpha$ and $\beta$ and this
solution is uniformly exponentially stable. By
Theorems~\ref{tConditionalPeriodic}
and~\ref{tConditionalStability}, the same is true for the original
problem~\eqref{eq2.1}, \eqref{eq2.3}.
\end{remark}

Assume that Condition~\ref{condmN0} holds. Let $z(x,t)$ be a
$2s$-periodic solution  of problem~\eqref{eq2.1}, \eqref{eq2.3}.
  Denote by
$\tilde\bz(t)=(z_0(t),\bz(t))$ the corresponding $2s$-periodic
solution
 of the guiding system~\eqref{eqInvariantDSN}. Let
us study the map $\tilde\Pi_\alpha$ and   the Poincar\'e map
$\tilde\Pi$ (see Sec.~\ref{secPeriodic}) of the guiding
system~\eqref{eqInvariantDSN}  in a neighborhood of
$\tilde\bz(0)$.

First of all, we consider  the projections of these operators onto
the $N$-dimensional space $V$ (see~\eqref{eqHDecomp}).

We consider  the orthogonal projector $$\bE:\tilde V\to V $$ given
by $\bE\tilde\bphi=\bphi$, where
$$
\tilde\bphi=\{\phi_j\}_{j=0}^N,\qquad \bphi =\{\phi_j\}_{j=1}^N.
$$

\begin{figure}[ht]
      {\ \hfill\epsfxsize65mm\epsfbox{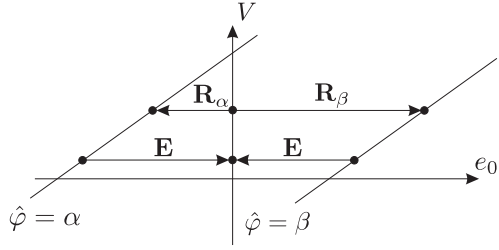}\hfill\ }
        \caption{The projection operator  $\bE$ and the lifting
operators $\bR_\alpha$ and $\bR_\beta$}
        \label{figProjER}
\end{figure}

We also introduce the lifting operator
$$\bR_\alpha:V \to \tilde V$$ given by
$$
\bR_\alpha(\bphi)=\left(\dfrac{\alpha}{m_0}-\dfrac{1}{m_0}\sum\limits_{k=1}^N
m_k\phi_k, \{\phi_j\}_{j=1}^N\right).
$$
Thus, $\bR_\alpha\bE(\tilde\bphi)=\tilde\bphi$ for
$\tilde\bphi\in\tilde V$ such that $\sum\limits_{j=0}^N m_j
\phi_j=\alpha$, and $\bE\bR_\alpha(\bphi)=\bphi$ for $\bphi\in V $
(see Fig.~\ref{figProjER}).

Denote by $\Pi_\alpha:V  \to V  $ the ``projection'' of
$\tilde\Pi_\alpha$ onto $V$ given by
$$
\Pi_\alpha(\bphi)=\bE\tilde\Pi_\alpha\bR_\alpha(\bphi).
$$
Similarly, one can define the operators $\bR_\beta$ and
$\Pi_\beta$.

The operators $\bE$, $\bR_\alpha$, and $\bR_\beta$ are
continuously (and even infinitely)  differentiable. Therefore, the
operators~$\Pi_\alpha$ and~$\Pi_\beta$ are also continuously
differentiable, provided so are $\tilde\Pi_\alpha$ and
$\tilde\Pi_\beta$.

\begin{figure}[ht]
      {\ \hfill\epsfxsize75mm\epsfbox{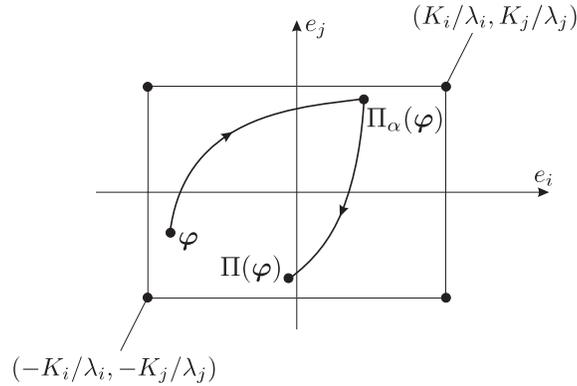}\hfill\ }
        \caption{The operators $\Pi_\alpha$ and $\Pi=\Pi_\beta\Pi_\alpha$ in the space
        $V =\Span(e_1,e_2,\dots,e_N)$}
        \label{figPiAlphaBeta}
\end{figure}

We introduce the operator $\Pi: V  \to V  $ by the formula
$$
\Pi(\bphi)=\bE\tilde\Pi\bR_\alpha(\bphi).
$$
 The following property of $\Pi$
is straightforward (see Fig.~\ref{figPiAlphaBeta}):
$$
\Pi=\Pi_\beta\Pi_\alpha.
$$

It is easy to see that the point $\bpsi=\bz(0)$ is a fixed point
of the map $\Pi$ acting in the $N$-dimensional space~$V$.

In the formulation of the following results, we will use the
following functions:
\begin{equation}\label{eqQs}
Q_j=Q_j(s)=\dfrac{2e^{-\lambda_j s}}{1+e^{-\lambda_j s}},\qquad
Q=Q(s)= m_0K_0+\sum\limits_{j=1}^Nm_jK_jQ_j(s).
\end{equation}

We note that, due to~\eqref{eqdvts}, we have at the switching
moment $s$
\begin{equation}\label{eqDerPialphas1}
\dfrac{d\hat z(t)}{dt}\Big|_{t=s}=Q(s).
\end{equation}
In particular, this implies that $Q(s)\ge 0$.

\begin{theorem}\label{tStablezStablePi}
Let Condition~$\ref{condmN0}$ hold, and let $z(x,t)$ be a
$2s$-periodic solution  of problem~$\eqref{eq2.1}$,
$\eqref{eq2.3}$. Assume that $Q(s)>0$. Then $z(x,t)$ is stable
{\rm (}uniformly exponentially stable{\rm )} if and only if the
fixed point $\bz(0)$ of the map $\Pi$ is so.
\end{theorem}
\proof Due to~\eqref{eqDerPialphas1}, we have $ \dfrac{d\hat
z(t)}{dt}\Big|_{t=s}>0 $. By symmetry, $ \dfrac{d\hat
z(t)}{dt}\Big|_{t=2s}<0 $. Now it remains to apply the formula
$\tilde\Pi^i(\tilde\bphi)=(\bR_\alpha\Pi^i\bE)(\tilde\bphi)$ and
Theorem~\ref{tConditionalStability}.
\endproof

To study the stability of the point $\bpsi=\bz(0)$, we consider
the  derivative of $\Pi$ at the point $\bpsi$.

\begin{lemma}\label{lDerPialphas}
Let Condition~$\ref{condmN0}$ hold. If $Q(s)>0$, then the operator
$\Pi_\alpha:V \to V $ is differentiable in a neighborhood of
$\bpsi=\bz(0)$ and the derivative
$$
D_\bpsi\Pi_\alpha(\bpsi):V\to V
$$
at the point $\bpsi=\bz(0)$ is given by
\begin{equation}\label{eqDerPiAlphas}
D_\bpsi\Pi_\alpha(\bpsi)\bphi=\sum\limits_{j=1}^N e^{-\lambda_j
s}\phi_je_j(x) +\dfrac{1}{Q(s)}\left(\sum\limits_{k=1}^N m_k
\left(1-e^{-\lambda_k s}\right)\phi_k\right) \sum\limits_{j=1}^N
K_jQ_j(s)e_j(x),
\end{equation}
where $e_1(x),\dots,e_N(x)$ form the basis in $V$ and $Q_j(s)$ and
$Q(s)$ are defined in~\eqref{eqQs}.
\end{lemma}
\proof Since $Q(s)>0$, it follows from~\eqref{eqDerPialphas1} that
$ \dfrac{d\hat z(t)}{dt}\Big|_{t=s}>0$. Therefore, applying
Lemma~4.2 in~\cite{GurHystUncond}, we have
$$
D_\bpsi\Pi_\alpha(\bpsi)\bphi=\sum\limits_{j=1}^N e^{-\lambda_j
s}\phi_je_j(x) +\left(\dfrac{d\hat
z(t)}{dt}\Big|_{t=s}\right)^{-1}\left(\sum\limits_{k=1}^N
m_k\left(1 - e^{-\lambda_k s}\right)\phi_k\right)
\sum\limits_{j=1}^N \lambda_j e^{-\lambda_j
s}\left(\dfrac{K_j}{\lambda_j}-\psi_j\right)e_j(x).
$$
Taking into account equalities~\eqref{eqphi_js}, \eqref{eqQs},
and~\eqref{eqDerPialphas1}, we obtain the desired
representation~\eqref{eqDerPiAlphas}.
\endproof

\begin{remark}\label{remA}
Due to Lemma~\ref{lDerPialphas}, the linear operator
$D_\bpsi\Pi_\alpha(\bpsi)$ is represented in the basis
$e_1(x),\dots,e_N(x)$ by the $(N\times N)$-matrix $\bA=\bA(s)$ of
the form
\begin{equation}\label{eqMatrA}
\bA=\left(
\begin{matrix}
1-E_1+S_1 \sigma_1 & S_1 \sigma_2       &  S_1 \sigma_3      & \dots &  S_1 \sigma_N\\
S_2 \sigma_1       & 1-E_2+S_2 \sigma_2 &  S_2 \sigma_3      & \dots &  S_2 \sigma_N\\
S_3 \sigma_1       & S_3 \sigma_2       & 1-E_3+S_3 \sigma_3 & \dots &  S_3 \sigma_N\\
\dots       & \dots      & \dots & \dots &  \dots\\
S_N \sigma_1       & S_N \sigma_2       &  S_N \sigma_3 & \dots &
1-E_N+S_N \sigma_N
\end{matrix}
\right),
\end{equation}
where
\begin{equation}\label{eqEj}
E_j=E_j(s)=1-e^{-\lambda_js},\qquad
S_j=S_j(s)=\dfrac{K_jQ_j(s)}{Q(s)},\qquad
\sigma_j=\sigma_j(s)=m_jE_j(s).
\end{equation}

Note that $\bA(0)$ is the identity matrix.
\end{remark}

The following lemma results from Lemma~\ref{lDerPialphas} and from
the symmetry of the periodic solution $z(x,t)$.

\begin{lemma}\label{lDerPis}
Let Condition~$\ref{condmN0}$ hold, and let $Q(s)>0$. Then the
operator $\Pi:V \to V $ is differentiable in a neighborhood of
$\bpsi=\bz(0)$ and the derivative
$$
D_\bpsi\Pi(\bpsi):V\to V
$$
at the point $\bpsi=\bz(0)$ is given in the basis
$e_1(x),\dots,e_N(x)$ by the matrix $\bA^2$, where $\bA$ is
defined in~\eqref{eqMatrA}.
\end{lemma}

Denote the eigenvalues of the matrix $\bA=\bA(s)$ by
$\mu_i=\mu_i(s)$, $i=1,\dots,N$.

The main result of this section is the following theorem. In
particular, we will use it to construct unstable periodic
solutions.

\begin{theorem}\label{tMuStability}
Let Condition~$\ref{condmN0}$ hold, and let $z(x,t)$ be a
 $2s$-periodic solution  of problem~$\eqref{eq2.1}$,
$\eqref{eq2.3}$. Assume that $Q(s)>0$.  Then the following
assertions are true.
\begin{enumerate}
\item All the eigenvalues $\mu_i$ of the matrix $\bA$ satisfy
$\mu_i\ne 1$. \item  If   $|\mu_i|<1$ for all $i=1,\dots N$, then
the $2s$-periodic solution $z(x,t)$ of problem~$\eqref{eq2.1}$,
$\eqref{eq2.3}$ is uniformly exponentially stable. \item If there
is an eigenvalue $\mu_k$ such that $|\mu_k|>1$, then the
$2s$-periodic solution $z(x,t)$ of problem~$\eqref{eq2.1}$,
$\eqref{eq2.3}$ is unstable.
\end{enumerate}
\end{theorem}
\proof Assertion 1 follows from Lemma~\ref{lDetA-I} below.  It is
known~\cite{Katok} that, under assumptions of items 2 and 3, a
 fixed point is, respectively, stable or unstable. By
Theorem~\ref{tStablezStablePi} and Lemma~\ref{lDerPis}, this fact
implies assertions 2 and 3.
\endproof

\begin{corollary}
Let Conditions~$\ref{condmKPositive}$ and~$\ref{condmN0}$ hold.
Then, for all sufficiently small $\beta-\alpha>0$, there exists a
  $2s$-periodic solution $z(x,t)$ of
problem~$\eqref{eq2.1}$, $\eqref{eq2.3}$ and assertions~$1$--$3$
in Theorem~$\ref{tMuStability}$ are true.
\end{corollary}
\proof The existence of $z(x,t)$ follows from
Theorem~\ref{tExistPeriodicSolution}. Moreover, we have shown in
the proof of Theorem~\ref{tExistPeriodicSolution} that $Q(s)>0$
for all sufficiently small $\beta-\alpha>0$, provided that
Condition~\ref{condmKPositive} holds. Thus, the hypothesis of
Theorem~\ref{tMuStability} are true. Therefore, the conclusions
are also true.
\endproof

Now we prove the following auxiliary result, which we have already
used in the proof of Theorem~\ref{tMuStability}.

\begin{lemma}\label{lDetA-I}
Let Condition~$\ref{condmN0}$ hold. If $Q(s)\ne 0$, then the
eigenvalues $\mu_i$ of $\bA$ satisfy
$$
\prod\limits_{i=1}^N(\mu_i-1)=(-1)^N
\dfrac{m_0K_0}{Q}\prod\limits_{i=1}^N E_i,
$$
where $Q$ is defined in~\eqref{eqQs} and $E_i$ in~\eqref{eqEj}.
\end{lemma}
\proof Substituting $\sigma_j=m_jE_j$, we have
$$
\prod\limits_{i=1}^N(\mu_i-1)=|\bA-\bI|=\prod\limits_{i=1}^N
E_i\cdot|\bB|,
$$
where
$$
\bB=\left(
\begin{matrix}
S_1 m_1-1 & S_1 m_2       &  S_1 m_3      & \dots &  S_1 m_N\\
S_2 m_1       & S_2 m_2-1 &  S_2 m_3      & \dots &  S_2 m_N\\
S_3 m_1       & S_3 m_2       & S_3 m_3-1 & \dots &  S_3 m_N\\
\dots       & \dots      & \dots & \dots &  \dots\\
S_N m_1       & S_N m_2       &  S_N m_3 & \dots & S_N m_N-1
\end{matrix}
\right)
$$
and $|\cdot|$ stands for the determinant of a matrix.

 Let us
compute the determinant of $\bB$:
$$
|\bB|=m_1\left|
\begin{matrix}
S_1   & S_1 m_2       &  S_1 m_3      & \dots &  S_1 m_N\\
S_2         & S_2 m_2-1 &  S_2 m_3      & \dots &  S_2 m_N\\
S_3         & S_3 m_2       & S_3 m_3-1 & \dots &  S_3 m_N\\
\dots       & \dots      & \dots & \dots &  \dots\\
S_N         & S_N m_2       &  S_N m_3 & \dots & S_N m_N-1
\end{matrix}\right|-
 \left|
\begin{matrix}
S_2 m_2-1 &  S_2 m_3      & \dots &  S_2 m_N\\
S_3 m_2       & S_3 m_3-1 & \dots &  S_3 m_N\\
\dots      & \dots & \dots &  \dots\\
S_N m_2       &  S_N m_3 & \dots & S_N m_N-1
\end{matrix}
\right|.
$$
To find the determinant of the first matrix, we multiply its first
column by $m_j$ and subtract it from the $j$th column for all
$j=2,\dots,N$. As a result, we have
$$
|\bB|=(-1)^{N-1}S_1m_1-
 \left|
\begin{matrix}
S_2 m_2-1 &  S_2 m_3      & \dots &  S_2 m_N\\
S_3 m_2       & S_3 m_3-1 & \dots &  S_3 m_N\\
\dots      & \dots & \dots &  \dots\\
S_N m_2       &  S_N m_3 & \dots & S_N m_N-1
\end{matrix}
\right|.
$$

Similarly decomposing the second determinant, we obtain (after
finitely many steps)
$$
|\bB|= (-1)^{N-1}(S_1m_1+\dots+S_Nm_N-1)=(-1)^N\dfrac{m_0K_0}{Q}.
$$
\endproof

\begin{remark}\label{remNinfinity}
Let us discuss modifications needed in the case of infinite  set
$\bbJ$ in Condition~\ref{condmN0}. The construction of the maps
$\Pi_\alpha,\Pi_\beta,\Pi$ is quite similar and the modifications
are obvious. The conclusion of Theorem~\ref{tStablezStablePi} with
the modified map $\Pi$ remains true.

Formula~\eqref{eqDerPiAlphas} for the Fr\'echet derivative
$D_\bpsi\Pi_\alpha(\bpsi)$ remains the same but the sums become
infinite. Their convergence follows from Remark~\ref{remAsymp}.
Formally, the linear operator $D_\bpsi\Pi_\alpha(\bpsi)$  can be
represented as the matrix~$\bA$ (see~\eqref{eqMatrA}), which now
becomes infinite-dimensional.

It is proved in~\cite{GurHystUncond} that the operators
$\Pi_\alpha,\Pi_\beta,\Pi$ are compact. Therefore, the same is
true for their Fr\'echet derivatives. In particular, this means
that the spectrum of $D_\bpsi\Pi_\alpha(\bpsi)$ consists of  no
more than countably many eigenvalues, which may accumulate only at
the origin. Thus, assertions 2 and 3 in Theorem~\ref{tMuStability}
remain true (possibly with $N=\infty$ in assertion 2).
\end{remark}

\subsection{Corollaries}

In this subsection, we assume that Condition~\ref{condmKPositive}
holds and that $\beta-\alpha>0$ and $s>0$ are sufficiently small.
Then $Q(s)>0$ and a $2s$-periodic solution $z(x,t)$ exists. Using
Theorem~\ref{tMuStability}, we provide some explicit conditions of
its stability or instability. Moreover, we will show that a
periodic solution may have a saddle structure.

The case $N=0$ is trivial (see Remark~\ref{remN=0}), so we begin
with the case $N=1$.

\begin{corollary}
Let Condition~$\ref{condmN0}$ hold with $N=1$. Then, for any
$\beta-\alpha>0$, there exists a unique periodic solution $z(x,t)$
of problem~\eqref{eq2.1}, \eqref{eq2.3}. The solution $z(x,t)$ is
uniformly exponentially stable.
\end{corollary}
\proof 1. By using the explicit formulas
(Rermark~\ref{rExplicitVj}) for the trajectories, we see that, for
any trajectory $v(x,t)$, the function $\hat v(t)$ either increases
for all $t>0$ or first decreases and than increases. In
particular, this implies that $dv/dt>0$ at the first switching
moment.

2. One can directly verify that the first characteristic function
$$
F(s):=m_0K_0s+2
m_1\dfrac{K_1}{\lambda_1}\cdot\dfrac{1-e^{-\lambda_1s}}{1+e^{-\lambda_1s}}
$$
satisfies one of the two conditions:
\begin{enumerate}
\item[(a)] $F(s)>0$ and increases for all $s>0$, or

\item[(b)] there is $s^*>0$ such that $F(s)<0$ for $0<s<s^*$ and
$F(s)>0$ and increases for all $s>s^*$.
\end{enumerate}
In both  cases, the equation $F(s)=\beta-\alpha$ has exactly one
positive root $s_1$.

3. Due to the observation in part 1 of the proof, the second
characteristic function
$$
H(t,s):=m_0k_0(t-s)+2
m_1\dfrac{K_1}{\lambda_1}\cdot\dfrac{e^{-\lambda_1s}-e^{-\lambda_1t}}{1+e^{-\lambda_1s}}=0
$$
 satisfies the inequality $H(t,s_1)<0$ for
all $t<s_1$. Therefore, by Theorem~\ref{tExistSymmetricSolution},
there is a unique $2s_1$ periodic solution of
problem~\eqref{eq2.1}, \eqref{eq2.3}.

3. To prove its stability, we note that the matrix $\bA$ consists
of one element $\mu_1$. It satisfies (due to Lemma~\ref{lDetA-I}
or by direct computation)
$$
\mu_1=1-E_1+S_1\sigma_1=1-(1-e^{-\lambda_1s_1})\dfrac{m_0K_0}{Q},
$$
where $Q>0$ due to~\eqref{eqDerPialphas1} and the observation in
part 1 of the proof. If we show that $\mu_1\in(-1,1)$, then the
stability result will follow from Theorem~\ref{tMuStability}.

Clearly, $\mu_1\ne 1$ for $s_1>0$. One can also show that
$\mu_1\ne -1$ for $s_1>0$. To do so, one can  check for example
that the equation $\mu_1=-1$ uniquely determines $m_1K_1$ as a
function of the other parameters. Then substituting it into the
formula for $F(s_1)$ yields the contradiction $F(s_1)<0$.

Since $\mu_1\ne\pm 1$, $\mu_1\in(-1,1)$ for sufficiently large
$s_1$, and $\mu_1$ continuously depends on $s_1$, it follows that
$\mu_1\in(-1,1)$ for any $s_1$.
\endproof

Now we consider the case $N=2$.
\begin{corollary}
Let Condition~$\ref{condmN0}$ hold with $N=2$, and let
\begin{equation}\label{eqmKN3}
M=m_0K_0+m_1K_1+m_2K_2>0.
\end{equation}
 Then, for all sufficiently small $\beta-\alpha>0$,
there exists a   $2s$-periodic solution $z(x,t)$ of
problem~\eqref{eq2.1}, \eqref{eq2.3} uniquely determined by
Theorem~$\ref{tExistPeriodicSolution}$. If
\begin{equation}\label{eqtrA2}
(M-m_1K_1)\lambda_1+ (M-m_2K_2)\lambda_2< 0,
\end{equation}
then $|\mu_1|,|\mu_2|>1$ and $z(x,t)$ is unstable for all
sufficiently small $\beta-\alpha>0$. If
\begin{equation}\label{eqtrA2'}
(M-m_1K_1)\lambda_1+ (M-m_2K_2)\lambda_2>0,
\end{equation}
then $|\mu_1|,|\mu_2|<1$ and $z(x,t)$ is exponentially stable for
all sufficiently small $\beta-\alpha>0$.
\end{corollary}
\proof 1. The matrix $\bA$ is a $(2\times2)$-matrix. Therefore, it
has two eigenvalues $\mu_1$ and $\mu_2$, which are either both
real or complex conjugate. Denote
$\delta=\delta_{1,2}=\mu_{1,2}-1$. Clearly, $\delta_{1,2}$ are the
eigenvalues of $\bA-\bI$; hence, they are the roots of the
quadratic equation
\begin{equation}\label{eqfordelta}
\delta^2-\tr(\bA-\bI)\delta+|\bA-\bI|=0.
\end{equation}

Let us compute $\tr(\bA-\bI)$ and $|\bA-\bI|$. Due
to~\eqref{eqMatrA} and~\eqref{eqEj},
$$
\tr(\bA-\bI) = E_1\left(-1 + \frac{m_1K_1Q_1}{Q}\right) + E_2\left(-1 + \frac{m_2K_2Q_2}{Q}\right).
$$
On the other hand, formulas \eqref{eqEj} and~\eqref{eqQs} imply
that $E_j=\lambda_js+O(s^2)$, and $Q_j(s)=1+O(s)$, and
$Q(s)=M+O(s)$. Therefore,
\begin{equation}\label{eqfortrA-I}
\tr(\bA-\bI)= s (\lambda_1 + O(s))\left(-1 +\frac{m_1K_1}{M} +
O(s)\right) + s (\lambda_2 + O(s))\left(-1 +\frac{m_2K_2}{M} +
O(s)\right) = -(Ls + O(s^2)),
\end{equation}
where
$$
L=M^{-1}((M-m_1K_1)\lambda_1+ (M-m_2K_2)\lambda_2).
$$

 Further, by Lemma~\ref{lDetA-I},
\begin{equation}\label{eqfordetA-I}
|\bA - \bI| =  E_1E_2 \frac{m_0K_0}{Q}=s^2(\lambda_1 +
O(s))(\lambda_2 + O(s))\left(\frac{m_0K_0}{M}+O(s)\right) = J^2s^2
+ O(s^3),
\end{equation}
where
$$
J^2=\lambda_1\lambda_2\dfrac{m_0K_0}{M}.
$$

It follows from~\eqref{eqfortrA-I} and \eqref{eqfordetA-I} that
Eq.~\eqref{eqfordelta} is equivalent to the following:
$$
\delta^2+(Ls+O(s^2))\delta+J^2s^2+O(s^3)=0.
$$
Thus,
$$
\mu_{1,2}=1-\dfrac{Ls}{2}\pm\dfrac{s\sqrt{L^2-4J^2+O(s)}}{2}+O(s^2).
$$

2. If  inequality~\eqref{eqtrA2} holds, then $L<0$ and
$\Re\mu_{1,2}>1$ for all small $s>0$.

Assume that inequality~\eqref{eqtrA2'} holds, i.e., $L>0$. If
$L^2-4J^2+O(s)\ge0$, then the eigenvalues $\mu_{1,2}$ are real and
belong to the interval $(0,1)$. If $L^2-4J^2+O(s)<0$, then
$\mu_{1,2}$ are complex conjugate and
$$
(\Re\mu_1)^2+(\Im\mu_1)^2=1-Ls+O(s^2)<1,
$$
i.e., $|\mu_{1,2}|<1$.
\endproof

\begin{example} Consider the problem described in Example~\ref{exClassification}.

Let  $m_0>0$, $m_1=m_2>0$, and $m_3=m_4=\dots=0$. Then
condition~\eqref{eqmKN3} holds. Therefore,
condition~\eqref{eqtrA2}, which implies the instability of the
 periodic solution for small $s$, takes the form
$$
\dfrac{m_0}{m_1}<\sqrt{2}\dfrac{\lambda_2-\lambda_1}{\lambda_2+\lambda_1}=\dfrac{3\sqrt{2}}{5},
$$
while condition~\eqref{eqtrA2'}, which implies the uniform
exponential stability of the   periodic solution for small $s$,
takes the form
$$
\dfrac{m_0}{m_1}>\sqrt{2}\dfrac{\lambda_2-\lambda_1}{\lambda_2+\lambda_1}=\dfrac{3\sqrt{2}}{5}.
$$
\end{example}

Finally, we show that periodic solutions can be unstable for $N\ge
3$. Moreover, if $N$ is odd, they may have a saddle structure.

\begin{corollary}
Let Condition~$\ref{condmN0}$ hold with  $N\ge 3$, and let
\begin{equation}\label{eqmKN}
M=\sum\limits_{j=0}^N m_jK_j>0.
\end{equation}
 Then, for all sufficiently small $\beta-\alpha>0$,
there exists a   $2s$-periodic solution $z(x,t)$ of
problem~\eqref{eq2.1}, \eqref{eq2.3} uniquely determined by
Theorem~$\ref{tExistPeriodicSolution}$. If
\begin{equation}\label{eqtrAN}
\sum\limits_{j=1}^N (M-m_jK_j)\lambda_j<0,
\end{equation}
then $z(x,t)$ is unstable.

If we additionally assume that $N$ is odd, then there is an
eigenvalue of $D_\bpsi\Pi(\bz(0))$ with real part greater than $1$
and a real eigenvalue in the interval $(0,1)$.
\end{corollary}
\proof 1. The matrix $\bA$ is an $(N\times N)$-matrix.  Due
to~\eqref{eqMatrA}, \eqref{eqEj}, and~\eqref{eqtrAN},
$$
\begin{aligned}
\sum\limits_{j=1}^N \mu_j=\tr\bA&=N+\sum\limits_{j=1}^N(S_jm_j-1)E_j\\
&=N-M^{-1}\sum\limits_{j=1}^N (M-m_jK_j)\lambda_js+O(s^2)>N
\end{aligned}
$$
for sufficiently small $s>0$. Therefore, the real part of at least
one eigenvalue is greater than $1$. By Theorem~\ref{tMuStability},
this implies the instability of $z(x,t)$.

2. Now we additionally assume that $N$ is odd. By
Lemma~\ref{lDetA-I},
$$
\prod_{j=1}^N(\mu_j-1)<0.
$$
Since $N$ is odd, the set of eigenvalues of $\bA$ consists of an
odd number of real eigenvalues $\mu_1,\dots,\mu_L$ ($1\le L\le N$)
and $(N-L)/2$ pairs of complex conjugate eigenvalues. Therefore,
$$
\prod_{j=1}^L(\mu_j-1)<0.
$$
Hence, there is at least one eigenvalue, e.g., $\mu_1$, which is
real and is less than $1$. Taking into account that $\mu_j(0)=1$
and $\mu_j(s)$ continuously depend on $s$, we see that $\mu_1\in
(0,1)$. Applying Lemma~\ref{lDerPis}, we complete the proof.
\endproof